\documentclass{amsart}
\usepackage{graphicx}
\usepackage{calc}
\usepackage{hyperref}

\usepackage[cmtip,arrow,all,2cell]{xy}
\usepackage{pb-diagram,pb-xy}
\usepackage{amsmath,amsthm,amssymb,tabularx}

\theoremstyle{plain}

\newtheorem{tm}{Theorem}[section]

\newtheorem{prop}{Proposition}[section]

\newtheorem{lm}{Lemma}[section]

\newtheorem{cor}{Corollary}[section]

\theoremstyle{definition}

\newtheorem{ex}{Example}[section]

\newtheorem{rem}[ex]{Remark}

\newcommand{\norm}[1]{\left\Vert#1\right\Vert}

\newcommand{\bit}{\begin{itemize}}
\newcommand{\eit}{\end{itemize}}
\newcommand{\btm}{\begin{tm}}
\newcommand{\etm}{\end{tm}}
\newcommand{\blm}{\begin{lm}}
\newcommand{\bprop}{\begin{prop}}
\newcommand{\eprop}{\end{prop}}
\newcommand{\elm}{\end{lm}}
\newcommand{\bcor}{\begin{cor}}
\newcommand{\ecor}{\end{cor}}
\newcommand{\bex}{\begin{ex}}
\newcommand{\eex}{\end{ex}}
\newcommand{\bcx}{\begin{cex}}
\newcommand{\ecx}{\end{cex}}
\newcommand{\bers}{\begin{ers}}
\newcommand{\eers}{\end{ers}}
\newcommand{\bdf}{\begin{df}}
\newcommand{\edf}{\end{df}}
\newcommand{\brem}{\begin{rem}}
\newcommand{\erem}{\end{rem}}
\newcommand{\bpr}{\begin{proof}}
\newcommand{\epr}{\end{proof}}

\def\le{\leqslant}
\def\ge{\geqslant}

\def \N {\mathbb{N}}
\def \Z {\mathbb{Z}}
\def \R {\mathbb{R}}

\def \C {\mathbb{C}}

\def\Epi{\operatorname{\sf Epi}}

\def\Ker{\operatorname{\sf Ker}}

\def\Mor{\operatorname{\sf Mor}}
\def\env{\operatorname{\sf env}}

\def\Env{\operatorname{\sf Env}}

\def\Spec{\operatorname{\sf{Spec}}}
\def\Ste{\operatorname{\sf{Ste}}}
\def\card{\operatorname{\sf{card}}}

\def\leftlim{\mathop{\varprojlim}\limits}
\def\rightlim{\mathop{\varinjlim}\limits}

\def \e {\varepsilon}
\def \ph {\varphi}

\def\Ker{\operatorname{\sf Ker}}
\def\Im{\operatorname{\sf Im}}
\def\im{\operatorname{\sf im}}
\def\Coim{\operatorname{\sf Coim}}
\def\coim{\operatorname{\sf coim}}
\def\red{\operatorname{\sf red}}
\def\id{\operatorname{\sf id}}
\def\Epi{\operatorname{\sf Epi}}
\def\DEpi{\operatorname{\sf DEpi}}

\def\sp{\operatorname{\sf span}}

\def\InvSteAlg{\operatorname{\sf InvSteAlg}}

\def\AugC*{\operatorname{\sf AugC^*}}

\begin{document}

\title{On continuous duality for Moore groups}

\author{S.S.Akbarov}

\maketitle

\begin{abstract}
In this paper we correct the errors of Yu.~N.~Kuznetsova's paper on the continuous duality for Moore groups.
\end{abstract}

In \cite{Kuznetsova} Yu.~N.~Kuznetsova made an attempt to construct a generalization of the Pontryagin duality to the class of all (not necessarily commutative) Moore groups. Its main result can be expressed in the diagram
 \begin{equation}\label{chetyrehugolnik-C-C*-0}
  \xymatrix @R=1.pc @C=2.pc
 {
 {\mathcal C}^\star(G)
 & \ar@{|->}[r]^{\Env} & &
 \Env{\mathcal C}^\star(G)
 \\
 & & &
 \ar@{|->}[d]^{\star}
 \\
 \ar@{|->}[u]^{\star}
 & & &
 \\
 {\mathcal C}(G)
 & &
 \ar@{|->}[l]_{\Env}
 &
 {\mathcal K}(G)
 },
 \end{equation}
where $G$ is an arbitrary Moore group, ${\mathcal C}(G)$ the algebra of continuous functions on $G$ (with the topology of uniform convergence on compact sets), ${\mathcal C}^\star(G)$ its dual stereotype space algebra considered as the group algebra of measures with compact support on $G$, $\Env$ the operation\footnote{For the operation of envelope we use a notation different from the one used in \cite{Kuznetsova}.} called in \cite{Kuznetsova} the {\it $C^*$-envelope}, which assigns to a topological algebra $A$ a new topological algebra $\Env A$ (see details in \cite{Kuznetsova}), $\Env{\mathcal C}^\star(G)$ the result of the application of this operation to the algebra ${\mathcal C}^\star(G)$, $\star$ the operation of taking the stereotype dual space, and ${\mathcal K}(G)$ the result of its application to the algebra $\Env{\mathcal C}^\star(G)$.

Later it turned out that the proof of this fact in \cite{Kuznetsova} contains some inaccuracies and errors. The author of this paper tried to correct them in \cite{Akbarov-C^infty-1,Akbarov-C^infty-2}, but as it was indicated in  \cite[Errata]{Akbarov-C^infty-2}, one of these errors was overlooked in \cite{Akbarov-C^infty-2}. As a corollary, by this moment the Kuznetsova theory can be  considered as proved (up to some specifications, in particular, a modification of the notion of envelope) only for a narrower class of the groups of the form $\R^n\times K\times D$, where $K$ is a compact group, and $D$ a discrete Moore group.

Here we correct the error that was made in \cite{Kuznetsova} (and repeated in \cite{Akbarov-C^infty-2}) with the help of the results of the author's paper \cite{Akbarov-ker-coker}.

The theory of envelopes of topological algebras with the applications in geometry was described by the author in \cite{Akbarov-stein-groups,Akbarov-env,Akbarov-C^infty-1,Akbarov-C^infty-2} (the results and the notations of the last three papers are used here). The special case of the Arens-Michael envelopes was considered in A.~Yu.~Pirkovaskii's papers \cite{Pirkovskii-DM,Pirkovskii-PAMS,Pirkovskii-MMO}.

\section{Preliminary results}

\subsection{SIN-groups and Moore groups.}

A set $U$ in a group $G$ is said to be {\it normal}, if it is invariant with respect to conjugations:
$$
a\cdot U\cdot a^{-1}\subseteq U,\qquad a\in G.
$$
A locally compact group $G$ is called a {\it SIN-group}\label{DEF:SIN-gruppa}, if its normal neighborhoods of unity form a local base. An equivalent definition: if the left and the right uniform structures on $G$ are equivalent. In particular, all these groups are unimodular.

The class of all SIN-groups contains the abelian groups, the compact groups and the discrete groups. The following proposition belongs to S.~Grosser and M.~Moskowitz \cite[2.13]{Grosser-Moskowitz-2}:

\btm\label{TH:stroenie-SIN}
Each SIN-group $G$ is a discrete extension of the group $\R^n\times K$, where $n\in\Z_+$, and $K$ is a compact group:
\begin{equation}\label{SIN-kak-rasshirenie}
1\to \R^n\times K=N\to G\to D\to 1
\end{equation}
(here $D$ is a discrete group).
\etm

A locally ccompact group $G$ is called a {\it Moore group}\label{DEF:Moore-gruppa}, if all its irreducible unitary (in the usual sense) representations are finite-dimensional.

\btm\label{TH:Moore->SIN}
Each Moore group is a SIN-group.\footnote{See \cite[p.1452]{Palmer}.}
\etm

\btm\label{TH:Moore->AM}
Eacj Moore group is amenable.\footnote{See \cite[p.1486]{Palmer}.}
\etm

\btm\label{TH:Moore->quotient}
Each (Hausdorff) quotient group $G/H$ of a Moore group $G$ is a Moore group.\footnote{This is obvious.}
\etm

\bcor\label{TH:G=Moore->D=Moore}
If $G$ is a Moore group, then in its representation \eqref{SIN-kak-rasshirenie} the group $D$ is also a Moore group (and, in particular, $D$ is amenable).
\ecor

\btm\label{TH:D=Moore->konech-rassh-Abelevoi}
Each discrete Moore group is a finite extension of an abelian group.\footnote{See  \cite[Theorem 12.4.26 and p.1397]{Palmer}.}
\etm

Let us call a locally compact group $G$ a {\it compact buildup of an abelian locally compact group}\label{DEF:komp-nadstr-abelevoi-gruppy}, if there exist closed subgroups $Z$ and $K$ in $G$ with the following properties:
\bit{

\item[1)] $Z$ is an abelian group,

\item[2)] $K$ is a compact group,

\item[3)] $Z$ and $K$ commute:
$$
\forall a\in Z,\quad \forall y\in K\qquad a\cdot y=y\cdot a,
$$

\item[4)] the product of $Z$ and $K$ is $G$:
$$
\forall x\in G\qquad \exists a\in Z\quad\exists y\in K\qquad x=a\cdot y.
$$
}\eit
If it is necessary to specify which groups in this constructions are used, then we say that $G$ is a {\it buildup of the abelian group $Z$ with the help of the compact group $K$}.

A {\it Lie-Moore group} is an arbitrary Lie group $G$, which is at the same time a Moore group.

\btm\footnote{Yu.~N.~Kuznetsova \cite{Kuznetsova}, see also \cite[Theorem 3.34]{Akbarov-C^infty-1}.}\label{TH:Lie-Moore}
Each Lie-Moore group $G$ is a finite extension of some compact buildup of an abelian group.
\etm

\btm\label{TH:Moore=lim-Lie-Moore}
Each Moore group $G$ is a projective limit of a system of Lie-Moore groups.
\etm
\bpr
By \cite[12.6.5]{Palmer}, $G$ is a projective limit of its quotient Lie groups $G/H$. By Theorem \ref{TH:Moore->quotient} they must be Moore groups.
\epr

\subsection{Hopf algebras in a monoidal category.}

We need several definitions from category theory. Suppose we have
 \bit{

\item[--] a partially ordered set $(I,\le)$,

\item[--] a category ${\tt K}$,

\item[--] a family $\{X^i;\ i\in I\}$ of objects in ${\tt K}$,

\item[--] a family $\{\iota_i^j:X^i\to X^j; \ i\le j\}$ of morphisms in ${\tt K}$,
 }\eit
and
 \bit{
\item[1)] the morphisms $\iota_i^i$ (where the upper and the lower indices coincide) are just identity morphisms:
$$
\iota_i^i=1_{X^i}
$$

\item[2)] for each indices $i\le j\le k$ the following diagram is commutative:
$$
\begin{diagram}
\node[2]{X^j}\arrow{se,t}{\iota_j^k}  \\
\node{X^i}\arrow[2]{e,b}{\iota_i^k}\arrow{ne,t}{\iota_i^j}
\node[2]{X^k}
\end{diagram}
$$
 }\eit\noindent
Then the family $\{X^i;\iota_i^j\}$ is called a {\it covariant system} in
the category ${\tt K}$ over the partially ordered set $(I,\le)$.

If $\{X_i,\iota_i^j\}$ and $\{Y_i,\varkappa_i^j\}$ are two covariant systems over a partially ordered set $(I,\le)$, then a {\it morphism} $\ph:\{X_i,\iota_i^j\}\to \{Y_i,\varkappa_i^j\}$ is a system of morphisms $\ph_i:X_i\to Y_i$ such that all the diagrams
$$
 \xymatrix 
 {
 X_i\ar[r]^{\iota_i^j}\ar[d]_{\ph_i}& X_j\ar[d]^{\ph_j} \\
 Y_i\ar[r]_{\varkappa_i^j} & Y_j
 }
$$
are commutative.

Let $\{X^i;\iota_i^j\}$ be a covariant system over a partially ordered set $I$ in a category ${\tt K}$.
 \bit{

\item[---] For each object $X$ in ${\tt K}$ a {\it projective cone} of a covariant system $\{X^i;\iota_i^j\}$ with the vertex $X$ is a system of morphisms $\pi^i:X \to X^i$ such that for all indices $i\le j$ the following diagram is commutative:
$$
\begin{diagram}
\node[2]{X}\arrow{sw,t}{\pi^i}\arrow{se,t}{\pi^j}  \\
\node{X^i}\arrow[2]{e,b}{\iota_i^j} \node[2]{X^j}
\end{diagram}
$$

\item[---] A projective cone $\{X,\pi^i\}$ of a covariant system $\{X^i;\iota_i^j\}$ is called a {\it projective limit}\label{DEF:proj-limit} of this system, if for each its other projective cone $\{Y,\rho^i\}$ there is a unique morphism $\tau:Y\to X$ such that for any index $i$ the following diagram is commutative:
\begin{equation}\label{DEF:proj-limit-diagr}
\begin{diagram}
\node{Y}\arrow{se,b}{\rho^i}\arrow[2]{e,t,--}{\tau}\node[2]{X}\arrow{sw,b}{\pi^i}
\\
\node[2]{X^i}
\end{diagram}
\end{equation}
In this case for the object $X$ and the morphisms $\pi^i$ and $\tau$ we use the notations:
$$
X=\lim_{\infty\gets j}X^j,\qquad \pi^i=\lim_{\infty\gets j}\iota_i^j,\qquad \tau=\lim_{\infty\gets j}\rho^j.
$$
 }\eit
The injective limits are defined dually.

Let us say that in a monoidal category ${\sf M}$ {\it the tensor product $\otimes$ commutes with the projective limits}, if the following natural identities hold:
$$
\lim_{\infty\gets (i,j)}(X_i\otimes Y_j)=(\lim_{\infty\gets i}X_i)\otimes (\lim_{\infty\gets j}Y_j).
$$

\bex\label{EX:odot-perestanov-s-proj-lim}
In the category $({\sf Ste},\odot)$ the tensor product $\odot$ commutes with the projective limits \cite[(2.4.39)]{Akbarov-C^infty-1}.
\eex

Similarly we say that in a monoidal category ${\sf M}$ {\it the tensor product $\otimes$ commutes with the injective limits}, if
$$
\lim_{(i,j)\to\infty}(X_i\otimes Y_j)=(\lim_{i\to\infty}X_i)\otimes (\lim_{j\to\infty}Y_j).
$$

\bex\label{EX:circledast-perestanov-s-inj-lim}
In the category $({\sf Ste},\circledast)$ the tensor product $\circledast$ commutes with injective limits \cite[(2.4.39)]{Akbarov-C^infty-1}.
\eex

\btm\label{TH:nasledovanie-Hopfa-predelami}
Suppose that in a monoidal category ${\sf M}$ the tensor product $\otimes$ commutes with the projective (respectively, injective) limits, and let $\{H_i;\iota_i^j\}$ be a covariant system in the category ${\sf Hopf}_{\sf M}$ of Hopf algebras in ${\sf M}$ over a directed in descending (respectively, in ascending) order set $I$. Then if the system  $\{H_i;\iota_i^j\}$ has a projective (respectively, injective) limit in the category ${\sf M}$, then this limit has a structure of Hopf algebra in ${\sf M}$ which turns it into a projective (respectively, injective) limit of the system $\{H_i;\iota_i^j\}$ in the category ${\sf Hopf}_{\sf M}$:
$$
\overset{{\sf M}}{\lim_{\infty\gets i}}H_i=\overset{{\sf Hopf}_{\sf M}}{\lim_{\infty\gets i}}H_i \qquad (\overset{{\sf M}}{\lim_{i\to\infty}}H_i=\overset{{\sf Hopf}_{\sf M}}{\lim_{i\to\infty}}H_i.
)
$$
\etm

\brem
Theorem \ref{TH:nasledovanie-Hopfa-predelami} remains true if we replace the Hopf algebras by bialgebras, or algebras, or coalgebras. We formulate it for Hopf algebras since in this form this proposition is used further in Corollary \ref{COR:proj-lim-Hopf-algebr-odot}, and later in Theorem \ref{TH:Env_C-C^*(G)-odot-Hopf}.
\erem

For proof we need the following two lemmas.

\blm\label{ph:lim-X_i->lim-Y_i}
Let $I$ be a partially ordered set, and $\{X_i,\iota_i^j\}$ and $\{Y_i,\varkappa_i^j\}$ two covariant systems over it with projective limits
$$
X=\lim_{\infty\gets i}X_i,\qquad Y=\lim_{\infty\gets i}Y_i.
$$
Then for each morphism of these covariant systems $\ph:\{X_i,\iota_i^j\}\to \{Y_i,\varkappa_i^j\}$ there exists a morphism between their projective limits $\ph:X\to Y$ such that all the following diagrams are commutative
$$
 \xymatrix 
 {
  & X\ar[dl]_{\pi_i}\ar[dr]^{\pi_j}\ar@{--}[d] & \\
 X_i\ar[rr]_(.7){\iota_i^j}\ar[dd]_{\ph_i}& \ar@{-->}[d]_{\ph} & X_j\ar[dd]^{\ph_j} \\
  & Y\ar[dl]_{\rho_i}\ar[dr]_{\rho_j} & \\
 Y_i\ar[rr]_{\varkappa_i^j} & & Y_j
 }
$$
\elm
\bpr
We have to notice that the system of morphisms $\{\ph_i\circ\pi_i\}$ is a projective cone for the covariant system $\{Y_i,\varkappa_i^j\}$, hence there exists a morphism $\ph$ such that all the remote lateral faces of the prism are commutatuve:
$$
\ph_i\circ\pi_i=\rho_i\circ\ph.
$$
The all other faces are commutative due to the definitions of the cone and of the morphism of the covariant systems.
\epr

Recall that a partially ordered set $(I,\le)$ is said to be {\it directed in descending order} (respectively, {\it in ascending order}), if for each $i,j\in I$ there is $k\in I$ such that
$$
k\le i\ \&\ k\le j
$$
(repectively, $i\le k\ \&\ j\le k$).

\blm\label{LM:lim-X_(i,i)=lim-X_(i,j)}
If $I$ is a set directed in descending order, and $\{X_{i,j}\}$ a covariant system over the cartesian square $I\times I$, having a projective limit, then the diagonal of this system $\{X_{i,i}\}$ has a projective limit as well, and these limits coincide:
\begin{equation}\label{lim-X_(i,i)=lim-X_(i,j)}
\lim_{\infty\gets i}X_{i,i}=\lim_{\infty\gets (i,j)}X_{i,j}.
\end{equation}
\elm
\bpr
This follows from the definition of the directed set.
\epr

\bpr[Proof of Theorem \ref{TH:nasledovanie-Hopfa-predelami}.] We have to define structure of Hopf algebra on the projective limit
$$
H=\overset{{\sf M}}{\lim_{\infty\gets i}}H_i.
$$
We can just use the formula which follows from \eqref{lim-X_(i,i)=lim-X_(i,j)}:
$$
\lim_{\infty\gets i}X_i\otimes X_i=\lim_{\infty\gets (i,j)}X_i\otimes X_j.
$$
This means for example that we can define the multiplication in $H$ just by passing from the multiplications $\mu_i:H_i\otimes H_i\to H_i$ to their projective limits by Lemma  \ref{ph:lim-X_i->lim-Y_i}
$$
\mu:\lim_{\infty\gets i}(H_i\otimes H_i)\to \lim_{\infty\gets i} H_i=H,
$$
then by Lemma \ref{LM:lim-X_(i,i)=lim-X_(i,j)} we replace $\lim_{\infty\gets i}(H_i\otimes H_i)$ by $\lim_{\infty\gets (i,j)}(H_i\otimes H_j)$, and after that we can use the commutativity of $\otimes$ with the projective limits, and replace $\lim_{\infty\gets (i,j)}(H_i\otimes H_j)$ by $\lim_{\infty\gets i}H_i\otimes\lim_{\infty\gets j}H_j$.
$$
\mu:H\otimes H=\lim_{\infty\gets i}H_i\otimes\lim_{\infty\gets j}H_j\cong\lim_{\infty\gets (i,j)}(H_i\otimes H_j)\cong\lim_{\infty\gets i}(H_i\otimes H_i)\to \lim_{\infty\gets i} H_i=H,
$$
The other structure morphisms are defined similarly.
\epr

Let us denote further by ${\sf Hopf}_\odot$ the class of Hopf algebras in the monoidal category $(\Ste,\odot,\C)$ of stereotype spaces, and by ${\sf Hopf}_\circledast$ the same class in the monoidal category $(\Ste,\circledast,\C)$. From the bicompleteness of the category {\sf Ste} \cite[Theorem 4.21]{Akbarov} and Examples \ref{EX:odot-perestanov-s-proj-lim} and \ref{EX:circledast-perestanov-s-inj-lim} we obtain the following important corollaries:

\bcor\label{COR:proj-lim-Hopf-algebr-odot}
Each covariant system $\{H_i;\iota_i^j\}$ in the category ${\sf Hopf}_\odot$ of stereotype Hopf algebras over a directed in descending order set $I$ has a projective limit, and this is a projective limit of the system $\{H_i;\iota_i^j\}$ in the category ${\sf Ste}$ with a proper structure of the Hopf algebra over $\odot$:
$$
\overset{{\sf Ste}}{\lim_{\infty\gets i}}H_i=\overset{{\sf Hopf}_\odot}{\lim_{\infty\gets i}}H_i.
$$
\ecor

\bcor\label{COR:inj-lim-Hopf-algebr-circledast}
Each covariant system $\{H_i;\iota_i^j\}$ in the category ${\sf Hopf}_\circledast$  of stereotype Hopf algebras over a directed in ascending order $I$ has an injective limit, and this is the injective limit of the system $\{H_i;\iota_i^j\}$ in the category ${\sf Ste}$ with a proper structure of projective Hopf algebra:
$$
\overset{{\sf Ste}}{\lim_{i\to\infty}}H_i=\overset{{\sf Hopf}_\circledast}{\lim_{i\to\infty}}H_i.
$$
\ecor

\subsection{Co-complete spaces.}
As usual, we say that a locally convex space $X$ is {\it complete}, if each Cauchy net in $X$ converges. The dual property in the theory of stereotype spaces is the following: a locally convex space $X$ is called {\it co-complete}\label{DEF:kopolno}, if each linear functional $f:X\to\C$ continuous on each totally bounded set $K\subseteq X$, is continuous on the whole $X$. The duality between these notions in the class of stereotype spaces is described by the equivalence \cite[Section 2.2]{Akbarov}:
$$
\text{$X$ is co-complete}\qquad\Longleftrightarrow\qquad \text{$X^\star$ is complete}.
$$

Further we use the following two definitions from \cite{Akbarov-env} ($\mathcal{U}(X)$ denotes the system of all neighborhoods of zero in $X$).
\bit{
\item
 A linear map of locally convex spaces $\ph:X\to Y$ is said to be {\it open}\label{DEF:open-map}, if the image $\ph(U)$ of each neighborhood of zero $U\subseteq X$ is a neighborhood of zero in the subspace $\ph(X)$ in $Y$ (with the topology induced from $Y$):
$$
\forall
U\in \mathcal{U}(X) \quad \exists V\in \mathcal{U}(Y) \quad \ph(U)\supseteq
\ph(X)\cap V.
$$

\item
A linear continuous map of locally convex spaces $\ph:X\to Y$ is said to be
{\it closed}, if for each totally bounded set $T\subseteq \overline{\ph(X)}\subseteq Y$ there is a totally bounded set $S\subseteq X$ such that $\ph(S)\supseteq T$. (Certainly, this implies in particularly, that the set of values $\ph(X)$ of the map $\ph$ is closed in $Y$.)
}\eit
In the class of stereotype spaces the openness and the closedness of a linear continuous map (i.e. of a morphism in this category) are dual properties \cite[Theorem 2.10]{Akbarov-env}:
$$
\text{$\ph:X\to Y$ is closed}\qquad\Longleftrightarrow\qquad \text{$\ph^\star:Y^\star\to X^\star$ is open}.
$$

\btm\label{TH:zamk-biektsija-v-kopolnoe-prostrancsvo}
Suppose $\ph:X\to Y$ is a closed bijective linear continuous map of stereotype spaces, and  $Y$ is co-complete. Then $\ph:X\to Y$ is an isomorphism of stereotype spaces.
\etm
\bpr
The space $X$ can be treated as a new, finer stereotype topologisation of the stereotype space $Y$, which preserves the system of totally bounded sets and the topology on each totally bounded set. If we pass to the dual spaces, then $Y^\star$ is a subspace in $X^\star$ with the topology induced from $X^\star$, and $Y^\star$ is dense in $X^\star$ (since the dual mapping $\ph:X\to Y$ is injective). At the same time $Y^\star$ is complete (since $Y$ is co-complete). Together these properties mean that $Y^\star$ and $X^\star$ coincide as locally convex spaces, hence the same is true for $Y$ and $X$.
\epr

\subsection{Lemma on epimorphism.}

Recall that in \cite{Akbarov-env} the notion of nodal decomposition of an arbitrary morphism $\ph:X\to Y$ was introduced. This is a representation of $\ph$ as a composition $$
\ph=\sigma\circ\beta\circ\pi,
$$
where $\sigma$ is a strong monomorphism, $\beta$ a bimorphism, and $\pi$ a strong epimorphism. If such a representation exists, then it is unique up to an isomorphism of its components, hence we can assign notations to the elements of this construction:
$$
\sigma=\im_\infty\ph,\quad \beta=\red_\infty\ph,\quad \pi=\coim_\infty\ph.
$$
The domain of the morphism $\sigma$ is denoted by $\Im_\infty\ph$, and the range of the morphism $\pi$ by $\Coim_\infty\ph$. Thus the morphism $\ph$ is decomposed as follows:
\begin{equation}\label{DEF:oboznacheniya-dlya-uzlov-razlozh}
\begin{diagram}
\node{X}\arrow{s,l}{\coim_\infty\ph}\arrow{e,t}{\ph}\node{Y} \\
\node{\Coim_\infty\ph}\arrow{e,t}{\red_\infty\ph}\node{\Im_\infty\ph}\arrow{n,r}{\im_\infty\ph}
\end{diagram}
\end{equation}
It is known \cite[Theorem 4.100]{Akbarov-env} that in the category $\Ste$ of stereotype spaces each morphism has nodal decomposition.

\blm\label{LM:epimorphism}
Let $\e:X\to Y$ and $\ph:Y\to Z$ be morphisms of stereotype spaces, and $\psi=\ph\circ\e$ be their composition. If $\e$ is an epimorphism, then there exists a unique morphism $\ph_\infty:Y\to\Im_\infty\psi$ such that the following diagram is commutative:
\begin{equation}\label{diagr:epimorphism}
 \xymatrix 
 {
 X\ar[dr]_{\psi_\infty}\ar[rr]^{\e}\ar@/_5ex/[ddr]_{\psi} & & Y\ar@{-->}[dl]^{\ph_\infty}
 \ar@/^5ex/[ddl]^{\ph} \\
  & \Im_\infty\psi \ar[d]_{\im_\infty\psi} & \\
   & Z  &
 }
\end{equation}
where $\psi_\infty=\red_\infty\psi\circ\coim_\infty\psi$.
\elm
\bpr
We need here formula \cite[(4.85)]{Akbarov-env}, which states that the nodal image of the map $\psi$ coincides with the envelope $\Env^Z\psi(X)$ (in the sense of \cite[$\S$ 2, (e)]{Akbarov-env}) of the image $\psi(X)$ of the mapping $\psi$:
$$
\im_\infty\psi=\Env^Z\psi(X).
$$
The morphism $\ph_\infty$ is built by the transfinite induction.

Let us describe in detail the zero step. Take a point $y\in Y$. Since $\e$ is an epimorphism, there is a net $\{x_i\}\subseteq X$ such that
$$
\e(x_i)\overset{Y}{\underset{i\to\infty}{\longrightarrow}}y.
$$
Hence
$$
\psi(x_i)=\ph(\e(x_i))\overset{Z}{\underset{i\to\infty}{\longrightarrow}}\ph(y).
$$
and therefore $\ph(y)\in\overline{\psi(X)}^Z$. This is true for each $y\in Y$, hence
$$
\ph(Y)\subseteq\overline{\psi(X)}^Z
$$
Now we can treat $\ph$ as a continuous mapping of the stereotype space $Y$ into the  (not necessarily stereotype, but pseudocomplete) space $\overline{\psi(X)}^Z$ (with the topology induced from $Z$):
$$
\ph: Y\to \overline{\psi(X)}^Z.
$$
Since the space $Y$ is pseudosaturated, from \cite[(1.26)]{Akbarov} we deduce that we can treat $\ph$ as a continuous mapping from $Y$ into the pseudosaturation $\Big(\overline{\psi(X)}^Z\Big)^\vartriangle$ of the space $\overline{\psi(X)}^Z$:
$$
\ph: Y\to\Big(\overline{\psi(X)}^Z\Big)^\vartriangle.
$$
Set $E_0=\Big(\overline{\psi(X)}^Z\Big)^\vartriangle$, and let $\ph_0$ and $\psi_0$ be the mappings $\ph$ and $\psi$ considered as having ranges in $E_0$. Then we obtain a diagram which is a zero approximation to \eqref{diagr:epimorphism}:
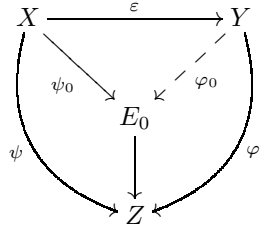
$$
 \xymatrix 
 {
 X\ar[dr]_{\psi_0}\ar[rr]^{\e}\ar@/_5ex/[ddr]_{\psi} & & Y\ar@{-->}[dl]^{\ph_0}
 \ar@/^5ex/[ddl]^{\ph} \\
  & E_0 \ar[d] & \\
   & Z  &
 }
$$

Further for each ordinal number $k$ we define $E_k$, $\ph_k$ and $\psi_k$ as follows:
 \bit{

\item[---] if $k$ is an isolated ordinal number, i.e. $k=j+1$ for some ordinal number $j$, then we apply the same trick as for $k=0$: we put  $E_k=\Big(\overline{\psi(X)}^{E_j}\Big)^\vartriangle$, and we obtain the diagram
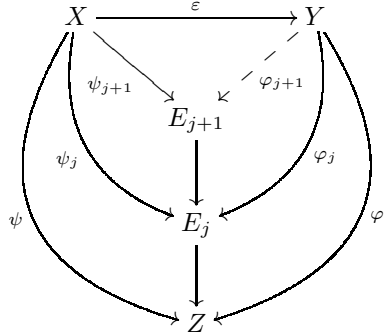
\begin{equation}\label{diagr:epimorphism-1}
 \xymatrix 
 {
 X\ar[dr]_{\psi_{j+1}}\ar[rr]^{\e}\ar@/_5ex/[ddr]_{\psi_j}\ar@/_10ex/[dddr]_{\psi} & & Y\ar@{-->}[dl]^{\ph_{j+1}}
 \ar@/^5ex/[ddl]^{\ph_j} \ar@/^10ex/[dddl]^{\ph} \\
  & E_{j+1} \ar[d] & \\
  & E_j \ar[d] & \\
   & Z  &
 }
\end{equation}

\item[---] if $k$ is a limit ordinal number, i.e. there is no such $j$ that $k=j+1$, then we put $E_k=\lim_{k\gets j} E_j$ (the projective limit in the category of stereotype spaces, and we obtain the diagram
$$
 
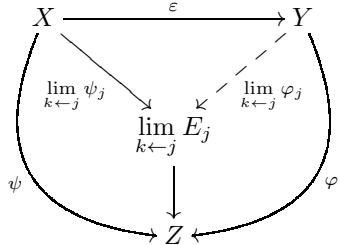
\xymatrix 
 {
 X\ar[dr]_{\lim\limits_{k\gets j}\psi_j}\ar[rr]^{\e}\ar@/_8ex/[ddr]_{\psi} & & Y\ar@{-->}[dl]^{\lim\limits_{k\gets j}\ph_j}
 \ar@/^8ex/[ddl]^{\ph} \\
  & \lim\limits_{k\gets j} E_j \ar[d] & \\
   & Z  &
 }
$$

}\eit
As a result we obtain a transfinite sequence of stereotype spaces $E_k$ and diagrams \eqref{diagr:epimorphism-1}. By \cite[(4.58)]{Akbarov-env}, this sequence stabilizes and its limit is exactly the diagram \eqref{diagr:epimorphism}.
\epr

\bcor\label{COR:Im_infty-ph-circ-e=Im_infty-ph}
Suppose $\e:X\to Y$ and $\ph:Y\to Z$ are morphisms of stereotype spaces, and $\e$ is an epimorphism. Then
\begin{equation}\label{Im_infty-ph-circ-e=Im_infty-ph}
\Im_\infty(\ph\circ\e)=\Im_\infty\ph
\end{equation}
\ecor
\bpr
On the one hand, we have an obvious implication
$$
(\ph\circ\e)(X)=\ph(\e(X))\subseteq \ph(X)\quad\Longrightarrow\quad
\Im_\infty(\ph\circ\e)\subseteq\Im_\infty\ph.
$$
On the other hand, Lemma \ref{LM:epimorphism} implies the embedding
$$
\ph(Y)\subseteq \Im_\infty(\ph\circ\e),
$$
which in its turn implies the embedding
$$
\Im_\infty\ph(Y)\subseteq \Im_\infty(\ph\circ\e).
$$
\epr

\section{Continuous envelope}

The notion of continuous envelope was introduced by the author in \cite{Akbarov-env} and was studied in detail in \cite{Akbarov-C^infty-2}. The exact definition is the following:

 \bit{
\item A {\it continuous envelope}\label{DEF:nepr-obolochka} $\env_{\mathcal C} A:A\to\Env_{\mathcal C} A$ of an involutive stereotype algebra $A$ is its envelope in the class $\DEpi$ of dense epimorphisms in the category $\InvSteAlg$ of involutive stereotype algebras with respect to the class $\Mor(\InvSteAlg,{\tt C}^*)$ of morphisms into $C^*$-algebras:
$$
\Env_{\mathcal C} A=\Env_{{\tt C}^*}^{\DEpi}A
$$
 }\eit

In detail, a {\it continuous extension} of an involutive stereotype algebra $A$ is a dense epimorphism $\sigma:A\to A'$ of involutive stereotype algebras such that for each  $C^*$-algebra $B$ and each involutive homomorphism $\ph:A\to B$ there is a (necessarily unique) homomorphism of stereotype algebras $\ph':A'\to B$ such that the following diagram is commutative:
\begin{equation}\label{DEF:diagr-nepr-rasshirenie}
 \xymatrix @R=2pc @C=1.2pc
 {
  A\ar[rr]^{\sigma}\ar[dr]_{\ph} & & A'\ar@{-->}[dl]^{\ph'} \\
  & B &
 }
\end{equation}
A {\it continuous envelope} of an involutive stereotype algebra $A$ is a continuous extension $\rho:A\to \Env_{\mathcal C} A$ such that for any other continuous extension $\sigma:A\to A'$ there is a (necessarily, unique) homomorphism of involutive stereotype algebras $\upsilon:A'\to \Env_{\mathcal C} A$ such that the following diagram is commutative:
$$
 \xymatrix @R=2pc @C=1.2pc
 {
  & A\ar[ld]_{\sigma}\ar[rd]^{\rho} &   \\
  A'\ar@{-->}[rr]_{\upsilon} &  & \Env_{\mathcal C} A
 }
$$

\subsection{Examples of continuous envelopes.}

Let us call a continuous mapping of topological spaces $\e:X\to Y$ a {\it covering}, if each compact set $T\subseteq Y$ is contained in the image of some compact set $S\subseteq X$. If the space $Y$ is Hausdorff, then this automatically implies that the mapping $\e$ is surjective. If in addition $\e$ is injective, then we call it an {\it exact covering}.\label{DEF:nalozhenie} In an exact covering $\e:X\to Y$ the space $Y$ can be treated as a new, coarser topologization of the space $X$, which does not change the system of compact sets and the topology on each compact set.

\btm\label{C-obolochka-podalgebry-v-C(M)} Let $A$ be an involutive stereotype subalgebra in the algebra $\mathcal{C}(M)$ of continuous functions on a paracompact locally compact space $M$, i.e. a (continuous and unital) monomorphism of involutive stereotype algebras is defined
$$
\iota:A\to \mathcal{C}(M).
$$
Then the continuous envelope of the algebra $A$ coincides with the algebra $\mathcal{C}(M)$
\begin{equation}\label{Env_C_A=C(M)}
\Env_{\mathcal C}  A=\mathcal{C}(M)
\end{equation}
(i.e. $\iota$ is a continuous envelope of $A$) if and only if the dual mapping of spectra $\iota^{\Spec}:\Spec(A)\gets M$ is an exact covering. \etm

The following fact was proved in \cite[Theorem 5.53]{Akbarov-env} (for the Kuznetsova envelopes in \cite[Theorem 2.11]{Kuznetsova}).

\btm\label{TH:Env_C-C^star(G)=C(widehat(G))}
The Fourier transform on a commutative locally compact group $H$
\begin{equation}\label{Fourier:C^star(H)->C(widehat(H))}
{\mathcal F}_H:{\mathcal C}^\star(H)\to{\mathcal C}(\widehat{H})\quad\Big|\quad
{\mathcal F}_H(\alpha)(\chi)=\alpha(\chi),
\end{equation}
is a continuous envelope of the group algebra ${\mathcal C}^\star(H)$. As a corollary,
\begin{equation}\label{Env_C-C^star(G)=C(widehat(G))}
\Env_{\mathcal C} {\mathcal C}^\star(H)={\mathcal C}(\widehat{H}).
\end{equation}
\etm

The following proposition is proved in \cite[Proposition 5.27]{Akbarov-C^infty-2}.

\btm\label{TH:env_C^star(Z-times-K)}
Let $Z$ be an abelian locally compact group, and $K$ a compact group. Then the formula
\begin{equation}\label{DEF:env_C^star(Z-times-K)}
(\Phi\delta^{(t,x)})_\sigma(\chi)=\chi(t)\cdot\sigma(x),\quad t\in Z,\ x\in K,\ \chi\in\widetilde{Z}, \ \sigma\in \widehat{K},
\end{equation}
defines a mapping
\begin{equation}\label{env_C^star(Z-times-K)}
\Phi:{\mathcal C}^\star(Z\times K)\to\prod_{\sigma\in\widehat{K}}{\mathcal C}\big(\widehat{Z},{\mathcal B}(X_\sigma)\big),
\end{equation}
which is a continuous envelope of the group algebra ${\mathcal C}^\star(Z\times K)$. As a corollary,
\begin{equation}\label{env_C-C^star(Z-times-K)}
\Env_{\mathcal C} {\mathcal C}^\star(Z\times K)={\mathcal C}\Big(\widehat{Z},\prod_{\sigma\in\widehat{K}}{\mathcal B}(X_\sigma)\Big)=
\prod_{\sigma\in\widehat{K}}{\mathcal C}\big(\widehat{Z},{\mathcal B}(X_\sigma)\big).
\end{equation}
\etm

The following theorem was proved in \cite[Theorem 2.7]{Akbarov-ker-coker}:

\btm\label{TH:Env_C-C*(Z-cdot-K)}
Let $Z\cdot K$ be a buildup of an abelian locally compact group $Z$ with the help of the compact group $K$. Then the continuous envelope $\Env_{\mathcal C} {\mathcal C}^\star(Z\cdot K)$ of its group algebra has the form
\begin{equation}\label{Env_C-C*(Z-cdot-K)}
\Env_{\mathcal C} {\mathcal C}^\star(Z\cdot K)\cong
\prod_{\sigma\in\widehat{K}}{\mathcal C}\Big(M_\sigma,{\mathcal B}(X_\sigma)\Big),
\end{equation}
where $\{M_\sigma;\ \sigma\in\widehat{K}\}$ is a family of closed subsets in the Pontryagin dual group $\widehat{Z}$ to the group $Z$.
\etm

\bprop\label{PROP:nepr-obolochka-C*-algebry}
If $A$ is a $C^*$-algebra, then $\Env_{\mathcal C} A=A$.
\eprop
\bpr
Certainly, the identity mapping $\id_A:A\to A$ is a continuous extension of the algebra $A$. If $\sigma:A\to A'$ is another continuous extension of the algebra $A$, then since  $A$ is a $C^*$-algebra, in the diagram
$$
\xymatrix @R=2.pc @C=3.0pc 
{
A\ar[rr]^{\sigma}\ar[dr]_{\id_A}& & A'\ar@{-->}[dl]^{\ph'} \\
& A &
}
$$
there is a unique dashed arrow $\ph'$. This means that $\id_A$ is a continuous envelope.
\epr

\bex\label{EX:Env_C-C_F=C_F}
If $F$ is a finite group, then it is convenient to denote its group algebra as $\C_F$ and to represent it as the dual space to the algebra $\C^F$ of all functions on $G$:
$$
\C_F=(\C^F)^\star.
$$
The group $F$ acts by shifts on the space of functions $L_2(F)$, hence the algebra $\C_F$ acts on $L_2(F)$ as well. This action can be understood as an embedding of $\C_F$ into the algebra of operators ${\mathcal B}(L_2(F))$. This algebra ${\mathcal B}(L_2(F))$ is a  $C^*$-algebra, hence $\C_F$ can also be treated as a $C^*$-algebra. As a corollary, by Proposition \ref{PROP:nepr-obolochka-C*-algebry}, its continuous envelope coincides with it:
\begin{equation}\label{Env_C-C_F=C_F}
\Env_{\mathcal C} \C_F=\C_F
\end{equation}
\eex

\subsection{Induced representation and the mapping $\Env_{\mathcal C}\theta:\Env_{\mathcal C} {\mathcal C}^\star(N)\to \Env_{\mathcal C} {\mathcal C}^\star(G)$.}

Recall the construction of the induced representation.
Let $N$ be an open normal subgroup in a locally compact group $G$. We put $F=G/N$ and consider $G$ as an extension of the group $N$ by the group $F$:
\begin{equation}\label{G-kak-rasshirenie}
\xymatrix 
{
1\ar[r] & N \ar[r]^{\eta} & G\ar[r]^{\ph} & F\ar[r]& 1
}
\end{equation}
(here  $\eta$ is the natural embedding, and $\ph$ the quotient mapping). Let us choose a function $\sigma:F\to G$ which is a coretraction for $\ph$,
\begin{equation}\label{ph(sigma(t))=t}
\ph(\sigma(t))=t,\qquad t\in F,
\end{equation}
and preserves the identity:
\begin{equation}\label{sigma(1_D)=1_G}
\sigma(1_F)=1_G.
\end{equation}
Then for each $g\in G$ the element $\sigma(\ph(g))$ belongs to the same coset of $N$ as $g$,
$$
g\in \sigma(\ph(g))\cdot N
$$
i.e.
\begin{equation}\label{g-sigma(ph(g))^(-1)-in-N}
g\cdot \sigma(\ph(g))^{-1}\in N,\qquad g\in G.
\end{equation}

Let us say that a mapping $\pi:G\to{\mathcal B}(X)$ is a {\it norm-continuous representation} of a group $G$ in a Hilbert space $X$ if it is continuous and satisfies the identities
\begin{equation}\label{pi(g-cdot-h)=pi(g)-cdot-pi(h)-nepr}
\pi(g\cdot h)=\pi(g)\cdot\pi(h),\qquad \pi(1_G)=1_B,\qquad g,h\in G,
\end{equation}

The following fact is a variant of Lemma 3.6 in \cite{Kuznetsova} (or of Theorem 3.38 in \cite{Akbarov-C^infty-2}).

\btm \label{TH:Ind-Preds}
Suppose a locally compact group $G$ is represented as an extension \eqref{G-kak-rasshirenie} of some open normal subgroup $N\subseteq G$, and   $\pi:N\to{\mathcal B}(X)$ is a norm-continuous representation. Consider the space $L_2(F,X)$ of square-summable functions $\xi:F\to X$ (with respect to the counting measure $\card$ on $F$). Then the formula
\begin{multline}\label{ind-representation}
\pi'(g)(\xi)(t)=\pi\Big(\underbrace{\sigma(t)\cdot g\cdot \sigma\big(\ph(\sigma(t)\cdot g)\big)^{-1}}_{\scriptsize\begin{matrix}\phantom{\quad\tiny \eqref{g-sigma(ph(g))^(-1)-in-N}}
\text{\rotatebox{90}{$\owns$}}{\quad\tiny \eqref{g-sigma(ph(g))^(-1)-in-N}}\\ N\end{matrix}}\Big)\Big(\xi\big(\underbrace{\ph(\sigma(t)\cdot g)}_{\scriptsize\begin{matrix}
\text{\rotatebox{90}{$\owns$}}\\ F\end{matrix}}\big)\Big)=\\=
\pi\Big(\sigma(t)\cdot g\cdot \sigma\big(t\cdot\ph(g))\Big)^{-1}\Big)\Big(\xi\big(t\cdot\ph(g))\Big),
 \qquad \xi\in L_2(F,X),\quad t\in F,\quad g\in G,
\end{multline}
defines a norm-continuous representation $\pi':G\to{\mathcal B}(L_2(F,X))$.
\etm

\bit{

\item The representation $\pi':G\to {\mathcal B}(L_2(F,X))$ defined in this way is called the {\it representation induced by the representation $\pi:N\to {\mathcal B}(X)$}.

}\eit

Let us call a {\it $C^*$-seminorm} on an involutive stereotype algebra $A$ a seminorm $p:A\to\R_+$ obtained as a composition of some (continuous, involutive and unital) homomorphism $\ph:A\to B$ into some $C^*$-algebra $B$ and the norm $\norm{\cdot}$ on $B$:
$$
p(x)=\norm{\ph(x)},\qquad x\in A.
$$
The set of all $C^*$-seminorms on $A$ is denoted by ${\mathcal P}(A)$. In a special case when $A={\mathcal C}^\star(G)$ we use the notation
$$
{\mathcal P}(G)={\mathcal P}\big({\mathcal C}^\star(G)\big).
$$

\bprop\label{PROP:otkrytost-Env-C^*(N)->Env-C^*(G)}
Each seminorm $p\in{\mathcal P}(N)$ is majorated by the restriction of some seminorm $q\in{\mathcal P}(G)$:
\begin{equation}\label{otkrytost-Env-C^*(N)->Env-C^*(G)}
p(\alpha)\le q(\alpha),\qquad\alpha\in {\mathcal C}^\star(N).
\end{equation}
\eprop
\bpr
We have to represent $p$ as a seminorm generated by some representation $\dot{\pi}:{\mathcal C}^\star(N)\to{\mathcal B}(X)$,
$$
p(\alpha)=\norm{\dot{\pi}(\alpha)},\qquad\alpha\in {\mathcal C}^\star(N),
$$
then consider the induced representation $\dot{\pi}':{\mathcal C}^\star(G)\to{\mathcal B}(L_2(D,X))$ and put
$$
q(\beta)=\norm{\dot{\pi}'(\beta)},\qquad\beta\in {\mathcal C}^\star(G).
$$
Then for a linear combination of delta-functions from $N$
\begin{equation}\label{alpha=sum_i-lambda_i-cdot-delta^(g_i)}
\alpha=\sum_i\lambda_i\cdot\delta^{g_i},\qquad g_i\in N,
\end{equation}
and for each $\xi\in L_2(F,X)$ we have
\begin{multline*}
\norm{\dot{\pi}'(\alpha)(\xi)(1)}= \norm{\dot{\pi}'(\sum_i\lambda_i\cdot\delta^{g_i})(\xi)(1)}=
\norm{\sum_i\lambda_i\cdot\dot{\pi}'(\delta^{g_i})(\xi)(1)}=
\norm{\sum_i\lambda_i\cdot\pi'(g_i)(\xi)(1)}=\\=
\norm{\sum_i\lambda_i\cdot\pi(\sigma(1)\cdot g_i\cdot \sigma(1)^{-1})(\xi(1))}=
\norm{\sum_i\lambda_i\cdot\pi(g_i)(\xi(1))}=
\norm{\sum_i\lambda_i\cdot\dot{\pi}(\delta^{g_i})(\xi(1))}=\\=
\norm{\dot{\pi}(\sum_i\lambda_i\cdot\delta^{g_i})(\xi(1))}=
\norm{\dot{\pi}(\alpha)(\xi(1))}
\end{multline*}
$$
\Downarrow
$$
$$
\norm{\dot{\pi}'(\alpha)(\xi)}=\sqrt{\sum_{t\in F}\norm{\dot{\pi}'(\alpha)(\xi)(t)}^2}\ge \norm{\dot{\pi}'(\alpha)(\xi)(1)}=\norm{\dot{\pi}(\alpha)(\xi(1))}
$$
$$
\Downarrow
$$
$$
q(\alpha)=\norm{\dot{\pi}'(\alpha)}=\sup_{\norm{\xi}\le 1}\norm{\dot{\pi}'(\alpha)(\xi)}\ge \sup_{\norm{\xi}\le 1}\norm{\dot{\pi}(\alpha)(\xi(1))}=
\sup_{\norm{\zeta}\le 1}\norm{\dot{\pi}(\alpha)(\zeta)}=
\norm{\dot{\pi}(\alpha)}=p(\alpha).
$$
The measures of the form \eqref{alpha=sum_i-lambda_i-cdot-delta^(g_i)} are dense in  ${\mathcal C}^\star(N)$, hence \eqref{otkrytost-Env-C^*(N)->Env-C^*(G)} holds for all $\alpha\in{\mathcal C}^\star(N)$.
\epr

Let further $\theta:{\mathcal C}^\star(N)\to {\mathcal C}^\star(G)$ denote the mapping induced by the embedding $\eta:N\subseteq G$ in \eqref{G-kak-rasshirenie}, and $\Env_{\mathcal C}\theta:\Env_{\mathcal C} {\mathcal C}^\star(N)\to \Env_{\mathcal C} {\mathcal C}^\star(G)$ the corresponding mapping of envelopes.
\begin{equation}\label{DIAGR:theta+Env-theta}
 \xymatrix  @R=2.pc @C=6.pc
{
 N\ar[r]^{\eta}\ar[d]_{\delta_N} & G \ar[d]^{\delta_{G}}
 \\
{\mathcal C}^\star(N)\ar[r]^{\theta}
 \ar[d]_{\env_{\mathcal C}} & {\mathcal C}^\star(G)
 \ar[d]^{\env_{\mathcal C}} \\
\Env_{\mathcal C}{\mathcal C}^\star(N)\ar[r]^{\Env_{\mathcal C}\theta} & \Env_{\mathcal C}{\mathcal C}^\star(G)
}
\end{equation}
Each seminorm $p\in{\mathcal P}(N)$ is uniquely extended to some seminorm on $\Env_{\mathcal C}{\mathcal C}^\star(N)$. We shall denote this extension by the same letter $p$. From Proposition \ref{PROP:otkrytost-Env-C^*(N)->Env-C^*(G)} we have

\bprop\label{PROP:otkrytost-Env-C^*(N)->Env-C^*(G)-1}
For each seminorm $p\in{\mathcal P}(N)$ there is a seminorm $q\in{\mathcal P}(G)$ such that
\begin{equation}\label{p(x)-le-q(Env-theta(x))}
p(x)\le q(\Env_{\mathcal C}\theta(x)),\qquad x\in\Env_{\mathcal C}{\mathcal C}^\star(N).
\end{equation}
\eprop
\bpr
Take $x=\env_{\mathcal C}\alpha$, where $\alpha\in {\mathcal C}^\star(N)$. Then for each seminorm $q$ from Proposition \ref{PROP:otkrytost-Env-C^*(N)->Env-C^*(G)} we have:
$$
p(x)=p(\env_{\mathcal C}\alpha)=p(\alpha)\le q\big(\theta(\alpha)\big)= q\Big(\env_{\mathcal C}\big(\theta(\alpha)\big)\Big)= q(\Env_{\mathcal C}\theta(\env_{\mathcal C}\alpha))=q(\Env_{\mathcal C}\theta(x))
$$
The elements of the form $x=\env_{\mathcal C}\alpha$, where $\alpha\in {\mathcal C}^\star(N)$, are dense in
$\Env_{\mathcal C}{\mathcal C}^\star(N)$, hence inequality \eqref{p(x)-le-q(Env-theta(x))} holds for all $x\in\Env_{\mathcal C}{\mathcal C}^\star(N)$.
\epr

\bprop\label{PROP:Env_C(theta)} Suppose in \eqref{G-kak-rasshirenie} $N=Z\cdot K$ is a compact buildup of an abelian locally compact group (and again an open normal subgroup in $G$). Then the continuous envelope $\Env_{\mathcal C}\theta:\Env_{\mathcal C} {\mathcal C}^\star(N)\to \Env_{\mathcal C} {\mathcal C}^\star(G)$ of the morphism $\theta:{\mathcal C}^\star(N)\to {\mathcal C}^\star(G)$ is an injective and an open\footnote{See definition at page \pageref{DEF:open-map}.} mapping of stereotype spaces.
\eprop
\bpr
We use here Corollary 2.2 from \cite{Akbarov-ker-coker}, according to which $\Env_{\mathcal C} {\mathcal C}^\star(N)$ is a locally convex projective limit of the quotient spaces ${\mathcal C}^\star(N)/p$, where $p\in{\mathcal P}(N)$:
\begin{equation}\label{E(C*(G))=LCS-leftlim}
\Env_{\mathcal C}{\mathcal C}^\star(N)={\tt LCS}\text{-}\kern-3pt\projlim_{p\in{\mathcal P}({\mathcal C}^\star(N))}{\mathcal C}^\star(N)/p
\end{equation}

1. From \eqref{E(C*(G))=LCS-leftlim} it follows that the seminorms $p\in{\mathcal P}(N)$ separate elements of $\Env_{\mathcal C} {\mathcal C}^\star(N)$. I.e. if $0\ne x\in \Env_{\mathcal C} {\mathcal C}^\star(N)$, then there is a seminorm $p\in{\mathcal P}(N)$ such that $p(x)>0$. By Proposition \ref{PROP:otkrytost-Env-C^*(N)->Env-C^*(G)} we can find a seminorm $q\in{\mathcal P}(G)$ such that \eqref{otkrytost-Env-C^*(N)->Env-C^*(G)} holds, and by \eqref{p(x)-le-q(Env-theta(x))} we have
$$
0<p(x)\le q(\Env_{\mathcal C}\theta(x)),
$$
hence $\Env_{\mathcal C}\theta(x)\ne 0$. This proves the injectivity of $\Env_{\mathcal C}\theta$.

2. From \eqref{E(C*(G))=LCS-leftlim} it follows also that the topology $\Env_{\mathcal C} {\mathcal C}^\star(N)$ is generated by seminorms $p\in{\mathcal P}(N)$ (without pseudosaturation). As a corollary we can think that the base neighborhoods of zero in $\Env_{\mathcal C} {\mathcal C}^\star(N)$ are generated by the seminorms $p\in{\mathcal P}(N)$. For every such a neighborhood of zero
$$
U=\{x\in \Env_{\mathcal C} {\mathcal C}^\star(N):\ p(x)\le\e \}
$$
by Proposition \ref{PROP:otkrytost-Env-C^*(N)->Env-C^*(G)} we choose a seminorm  $q\in{\mathcal P}(G)$ such that \eqref{p(x)-le-q(Env-theta(x))} holds, and then for a neighborhood of zero
$$
V=\{y\in \Env_{\mathcal C}{\mathcal C}^\star(G):\ q(y)\le\e \}
$$
we obtain
\begin{multline*}
y\in \Env_{\mathcal C}\theta\big(\Env_{\mathcal C}{\mathcal C}^\star(N)\big)\cap V\quad\Longrightarrow\quad
\exists x\in \Env_{\mathcal C}{\mathcal C}^\star(N)\quad y=\Env_{\mathcal C}\theta(x)\in V\quad\Longrightarrow \\ \Longrightarrow\quad p(x)\le q\big(\Env_{\mathcal C}\theta(x)\big)\le\e\quad\Longrightarrow\quad x\in U \quad\Longrightarrow\quad y=\Env_{\mathcal C}\theta(x)\in\Env_{\mathcal C}\theta(U).
\end{multline*}
This proves the openness of $\Env_{\mathcal C}\theta$.
\epr

\subsection{Algebra ${\mathcal K}(G)$}

Recall (see \cite{Akbarov-C^infty-2}) that on each locally compact group $G$ the algebra  ${\mathcal K}(G)$ is defined by the formula
 \begin{equation}\label{DEF:K(G)}
{\mathcal K}(G):=\Big(\Env_{\mathcal C} {\mathcal C}^\star(G)\Big)^\star.
 \end{equation}

To each element $u\in {\mathcal K}(G):=\Big(\Env_{\mathcal C} {\mathcal C}^\star(G)\Big)^\star$ one can assign a chain of mappings
$$
G\overset{\delta}{\longrightarrow}{\mathcal C}^\star(G)\overset{\env_{\mathcal C} {{\mathcal C}^\star(G)}}{\longrightarrow}\Env_{\mathcal C} {\mathcal C}^\star(G)\overset{u}{\longrightarrow}\C.
$$

\btm\label{TH:K(G)->C(G)} The mapping $u\mapsto u\circ \env_{\mathcal C} {{\mathcal C}^\star(G)}\circ\delta$ coincides with the mapping $\big(\env_{\mathcal C} {\mathcal C}^\star(G)\big)^\star$ dual to the mapping $\env_{\mathcal C} {{\mathcal C}^\star(G)}$:
\begin{equation}\label{K(G)->C(G)}
\big(\env_{\mathcal C}{\mathcal C}^\star(G)\big)^\star(u)=u\circ \env_{\mathcal C} {{\mathcal C}^\star(G)}\circ\delta
\end{equation}
and injectively and homomorphically embeds ${\mathcal K}(G)$ into ${\mathcal C}(G)$ as an involutive subalgebra (and as a corollary, the operations of summing, multiplication and involution in ${\mathcal K}(G)$ are pointwise).
\etm

The following proposition was proved in \cite[Theorem 5.40]{Akbarov-C^infty-2}.

\btm\label{LM:sdvig-v-K(G)}
The shift (the left and the right) by an arbitrary element $a\in G$ is an isomorphism of the stereotype algebra ${\mathcal K}(G)$.
\etm

\bit{
 \item[$\bullet$] An {\it involutive character} on an involutive stereotype algebra $A$ over $\C$ is an arbitrary (continuous, involutive and unital) homomorphism $s:A\to \C$. The space of all involutive characters on $A$ with the topology of uniform convergence on totally bounded sets in $A$ is called the {\it involutive spectrum}\label{DEF:Spec_R(G)} (or just {\it spectrum}) of $A$, and is denoted by $\Spec(A)$.
 }\eit

The following lemmas are used further in Theorem \ref{TH:Spec-K(G)=G}.

\blm\label{LM:Spec-K(G)-dlya-diskretnoi-gruppy}
If $G$ is an amenable discrete group, then the mapping of spectra $G\to \Spec{\mathcal K}(G)$ is a bijection.
\elm
\bpr This is Lemma 5.53 in \cite{Akbarov-C^infty-2}. \epr

\blm\label{LM:Spec(R^n-times-K)}
For each compact group $K$ and for each $n\in\N$ the mapping of spectra $\R^n\times K\to \Spec{\mathcal K}(\R^n\times K)$ is a homeomorphism.
\elm
\bpr
This is Lemma 5.46 in \cite{Akbarov-C^infty-2}.
\epr

\blm\label{LM:1_L-in-K(G)}
Let $N$ be an open normal subgroup in a locally compact group $G$. Then for each coset $L\in G/N$ its characteristic function $1_L$ is an element of the space ${\mathcal K}(G)$:
\begin{equation}\label{1_L-in-K(G)}
1_L\in {\mathcal K}(G)
\end{equation}
\elm
\bpr
This is proved similarly with Lemma 5.47 in \cite{Akbarov-C^infty-2}.
\epr

For each coset $L\in G/N$ we use the notations
\begin{equation}\label{DEF:K_L(G)}
{\mathcal K}_L(G)=1_L\cdot{\mathcal K}(G),\qquad {\mathcal K}_{G\setminus L}(G)=(1-1_L)\cdot{\mathcal K}(G)
\end{equation}
(here 1 is the identity in the algebra ${\mathcal K}(G)$). We endow these spaces with the structure of immediate subspaces in ${\mathcal K}(G)$ (see \cite{Akbarov-env}). From  \eqref{1_L-in-K(G)} we have

\blm\label{LM:K_L(G)-K_G-L(G)} Let $N$ be an open normal subgroup in a locally compact group $G$. Then the spaces ${\mathcal K}_L(G)$ and ${\mathcal K}_{G\setminus L}(G)$ complement each other in ${\mathcal K}(G)$:
\begin{equation}\label{K_N+K_(G-N)=K(G)}
{\mathcal K}_L(G)\oplus {\mathcal K}_{G\setminus L}(G)={\mathcal K}(G)
\end{equation}
(i.e. ${\mathcal K}(G)$ is a direct sum in the category of stereotype spaces).
\elm

Consider the morphism $\Env_{\mathcal C}\theta:\Env_{\mathcal C} ({\mathcal C}^\star(N))\to \Env_{\mathcal C} {\mathcal C}^\star(G)$ in diagram \eqref{DIAGR:theta+Env-theta} and denote by $\psi$  its dual map: $\psi=(\Env_{\mathcal C}\theta)^\star:{\mathcal K}(N)\gets {\mathcal K}(G)$.

\blm\label{LM:K(N)=K_N(G)} Suppose in \eqref{G-kak-rasshirenie} $N=Z\cdot K$ is a compact buildup of an abelian group (and again an open normal subgroup in the locally compact group $G$). Then the morphism of stereotype spaces $\psi=(\Env_{\mathcal C}\theta)^\star:{\mathcal K}(N)\gets {\mathcal K}(G)$ has the following properties:
\bit{
\item[(i)] its kernel is the second component in the decomposition \eqref{K_N+K_(G-N)=K(G)} (with $L=N$):
\begin{equation}\label{ph(K_(G-N)(G))=0}
\Ker\psi={\mathcal K}_{G\setminus N}(G).
\end{equation}
\item[(ii)] the restriction $\psi|_{{\mathcal K}_N(G)}:{\mathcal K}_N(G)\to {\mathcal K}(N)$ is an isomorphism of stereotype algebras.
}\eit
\elm
\bpr
1. For (i) let us consider the diagram
$$
\xymatrix @R=2.pc @C=8.0pc 
{
{\mathcal C}^\star(N)\ar[r]^{\env_{\mathcal C} {{\mathcal C}^\star(N)}}\ar[d]_{\theta} & \Env_{\mathcal C} ({\mathcal C}^\star(N))\ar[d]^{\Env_{\mathcal C}\theta}\\
{\mathcal C}^\star(G)\ar[r]^{\env_{\mathcal C} {{\mathcal C}^\star(G)}} & \Env_{\mathcal C} {\mathcal C}^\star(G)
}
$$
If $u\in \Ker\psi$, then we have a chain
\begin{multline*}
0=\psi(u)=u\circ \Env_{\mathcal C}\theta\quad\Longrightarrow\quad 0=u\circ \Env_{\mathcal C}\theta\circ \env_{\mathcal C} {{\mathcal C}^\star(N)}=
u\circ \env_{\mathcal C} {{\mathcal C}^\star(G)}\circ \theta \quad\Longrightarrow\\ \Longrightarrow\quad 0=u\circ \env_{\mathcal C} {{\mathcal C}^\star(G)}\circ \theta\circ\delta^N \quad\Longrightarrow\quad 0=u\Big|_N \quad\Longrightarrow\\ \Longrightarrow\quad 0=u\cdot 1_N
\quad\Longrightarrow\quad u=u\cdot (1-1_N) \quad\Longrightarrow\quad u\in {\mathcal K}_{G\setminus N}(G).
\end{multline*}
On the contrary, if $u\in {\mathcal K}_{G\setminus N}(G)$, then
\begin{multline*}
0=\psi(u)=u\circ \Env_{\mathcal C}\theta\quad\overset{\env_{\mathcal C} {{\mathcal C}^\star(N)}\in\Epi}{\Longleftarrow}\quad 0=u\circ \Env_{\mathcal C}\theta\circ \env_{\mathcal C} {{\mathcal C}^\star(N)}=
u\circ \env_{\mathcal C} {{\mathcal C}^\star(G)}\circ \theta \quad\overset{\overline{\sp\delta^N}={\mathcal C}^\star(N)}{\Longleftarrow}\\ \Longleftarrow\quad 0=u\circ \env_{\mathcal C} {{\mathcal C}^\star(G)}\circ \theta\circ\delta^N \quad\Longleftarrow\quad 0=u\Big|_N \quad\Longleftarrow\\ \Longleftarrow\quad 0=u\cdot 1_N
\quad\Longleftarrow\quad \exists v\in {\mathcal K}(G)\quad u=v\cdot (1-1_N) \quad\Longleftarrow\quad u\in {\mathcal K}_{G\setminus N}(G).
\end{multline*}

2. Let us prove (ii). Recall that by Proposition \ref{PROP:Env_C(theta)} the mapping  $\Env_{\mathcal C} \theta:\Env_{\mathcal C} {\mathcal C}^\star(N)\to \Env_{\mathcal C} {\mathcal C}^\star(G)$ is injective and open. This implies that the dual mapping $\psi:{\mathcal K}(G)\to {\mathcal K}(N)$ is closed, and hence, its coimage  $\psi|_{{\mathcal K}(G)}:{\mathcal K}(G)/\Ker\psi={\mathcal K}_N(G)\to {\mathcal K}(N)$ is again a closed mapping. On the other hand, since $\Env_{\mathcal C}\theta$ injective,  $\psi=(\Env_{\mathcal C}\theta)^\star$ is an epimorphism of stereotype spaces. Thus,  $\psi$ is a closed epimorphism. This implies that $\psi$ is surjective. Let us note in addition that $\psi$ is injective, since
$$
u\in {\mathcal K}_N(G)\ \&\ \psi(u)=0\quad\Longrightarrow\quad u\in u\in {\mathcal K}_N(G)\cap\Ker\psi={\mathcal K}_N(G)\cap{\mathcal K}_{G\setminus N}(G)=0\quad\Longrightarrow\quad u=0.
$$
Finally, one more important note: the space ${\mathcal K}(N)$ is co-complete (in the sense of definition at page \pageref{DEF:kopolno}), since its dual space
$$
 {\mathcal K}(N)^\star =\Env_{\mathcal C} {\mathcal C}^\star(N)=\eqref{Env_C-C*(Z-cdot-K)}=\prod_{\sigma\in\widehat{K}}{\mathcal C}\Big(M_\sigma,{\mathcal B}(X_\sigma)\Big),
$$
is complete. Thus, $\psi|_{{\mathcal K}_N(G)}:{\mathcal K}_N(G)\to {\mathcal K}(N)$ is a closed bijective continuous mapping of stereotype spaces, and its range is a co-complete space. By Theorem \ref{TH:zamk-biektsija-v-kopolnoe-prostrancsvo}, $\psi|_{{\mathcal K}_N(G)}:{\mathcal K}_N(G)\to {\mathcal K}(N)$ is an isomorphism of stereotype spaces.
\epr

The {\it involutive spectrum} $\Spec A$ on an involutive stereotype algebra $A$ is the set of all involutive characters, i.e. involutive (continuous and unital) homomorphisms  $\chi:A\to\C$, endowed with the topology of uniform convergence on compact sets in $A$.

\btm\label{TH:Spec-K(G)=G}
If $G$ is a Moore group, then the involutive spectrum of the algebra ${\mathcal K}(G)$ is topologically isomorphic to $G$:
\begin{equation}\label{Spec-K(G)=G}
\Spec{\mathcal K}(G)=G
\end{equation}
\etm
\bpr Let $G$ be a Moore group.

1. Let us first show that the mapping of spectra $G\to \Spec{\mathcal K}(G)$ is a surjection. Let $\chi:{\mathcal K}(G)\to\C$ be an involutive character. The homomorphism  $G\to D$ from \eqref{SIN-kak-rasshirenie} generates a homomorphism ${\mathcal C}^\star(G)\to {\mathcal C}^\star(D)$, which generates a homomorphism $\Env_{\mathcal C} {\mathcal C}^\star(G)\to \Env_{\mathcal C} {\mathcal C}^\star(D)$, and finally a homomorphism ${\mathcal K}(G)\gets {\mathcal K}(D)$. Let us denote it by $\ph:{\mathcal K}(D)\to {\mathcal K}(G)$. The composition $\chi\circ\ph:{\mathcal K}(D)\to \C$ is an involutive continuous character on ${\mathcal K}(D)$, and $D$ is a Moore group (by Corollary \ref{TH:G=Moore->D=Moore}), hence an amenable group (by Theorem \ref{TH:Moore->AM}). By Lemma \ref{LM:Spec-K(G)-dlya-diskretnoi-gruppy} $\chi$ is a delta-function:
$$
(\chi\circ\ph)(u)=u(L),\qquad u\in {\mathcal K}(D),
$$
for some $L\in G/N$. Consider the space ${\mathcal K}_L(G)$ from \eqref{DEF:K_L(G)} and denote by $\rho_L$ its embedding into ${\mathcal K}(G)$. Let us also denote by $\sigma$ the embedding ${\mathcal K}(N)\to {\mathcal K}_N(G)$, i.e. the isomorphism defined by Lemma \ref{LM:K(N)=K_N(G)}(ii). Take $b\in L$, then $L=N\cdot b$, and let $\tau_b:{\mathcal K}(G)\to{\mathcal K}(G)$ be the shift by the element $b^{-1}$ (acting on ${\mathcal K}(G)$ by Theorem \ref{LM:sdvig-v-K(G)}):
$$
\tau_bu=b^{-1}\cdot u,\qquad {\mathcal K}(G).
$$
It maps the space ${\mathcal K}(N)$ into the space ${\mathcal K}_L(G)$, hence a mapping is defined $\sigma_L=\tau_b\circ\sigma:{\mathcal K}(N)\to {\mathcal K}_L(G)$. Set
$$
\chi_L=\chi\circ\rho_L,\qquad \chi_N=\chi_L\circ\sigma_L
$$
and denote by $\rho_L:{\mathcal K}_L(G)\to {\mathcal K}(G)$ the natural embedding. We obtain a commutative diagram
$$
\xymatrix @R=2.pc @C=3.0pc 
{
{\mathcal K}(N)\ar[r]^{\sigma_L}\ar[dr]_{\chi_N} & {\mathcal K}_L(G)\ar[d]^{\chi_L}\ar[r]^{\rho_L} & {\mathcal K}(G)\ar[dl]^{\chi}\\
 & \C &
}
$$
Since $\chi_N$ is a character on ${\mathcal K}(N)$ by Lemma \ref{LM:Spec(R^n-times-K)} it must be a delta-function:
\begin{equation}\label{chi_N(u)=u(a)}
\chi_N(u)=u(a),\qquad u\in{\mathcal K}(N)
\end{equation}
for some $a\in N$. Then
\begin{multline*}
\chi(u)=\chi(1_L)\cdot\chi(u)=\chi(1_L\cdot u)=\chi_L(1_L\cdot u)=\chi_N(\sigma_L^{-1}(1_L\cdot u))=\eqref{chi_N(u)=u(a)}=
\sigma_L^{-1}(1_L\cdot u)(a)=\sigma(\sigma_L^{-1}(1_L\cdot u))(a)=\\=
(\sigma\circ\sigma_L^{-1})(1_L\cdot u))(a)=(\sigma\circ(\tau_b\circ\sigma)^{-1})(1_L\cdot u))(a)=
(\sigma\circ\sigma^{-1}\circ\tau_b^{-1})(1_L\cdot u))(a)=(\sigma\circ\sigma^{-1}\circ\tau_b^{-1})(1_L\cdot u))(a)=\\=
\tau_{b^{-1}}(1_L\cdot u))(a)=(b\cdot(1_L\cdot u))(a)=(1_L\cdot u)(\underbrace{a\cdot b}_{\scriptsize\begin{matrix}\text{\rotatebox{90}{$\owns$}}\\ L\end{matrix}})=u(a\cdot b)=\delta^{a\cdot b}(u).
\end{multline*}

2. Now let us verify that the mapping of spectra $G\to \Spec{\mathcal K}(G)$ is an injection. Take $a\ne b\in G$. If $a\cdot b^{-1}\notin N$, i.e. $a\notin b\cdot N$, then the characteristic function $1_L\in{\mathcal K}(G)$ of the class $L=b\cdot N$ from Lemma  \ref{LM:1_L-in-K(G)} distinguishes $a$ and $b$:
$$
1_L(a)=0\ne 1=1_L(b).
$$
Suppose $a\in b\cdot N$, i.e. $a\cdot b^{-1}\in N$. Then by Lemma  \ref{LM:Spec(R^n-times-K)} we can choose a function $u\in{\mathcal K}(N)$ such that
$$
u(a\cdot b^{-1})\ne u(1).
$$
By Lemma \ref{LM:K(N)=K_N(G)} there is a function $v\in{\mathcal K}_N(G)$ such that  $u\big|_N=v\big|_N$, and therefore
$$
v(a\cdot b^{-1})\ne v(1).
$$
By Theorem \ref{LM:sdvig-v-K(G)} the shift $b^{-1}\cdot v$ again lies in ${\mathcal K}(G)$, and for this function we have
$$
(b^{-1}\cdot v)(a)=v(a\cdot b^{-1})\ne v(1)=v(b\cdot b^{-1})=(b^{-1}\cdot v)(b).
$$

3. It remains to verify the openness of the mapping $G\to \Spec{\mathcal K}(G)$. Suppose that $a_i\to a$ in $\Spec{\mathcal K}(G)$. From Theorem \ref{LM:sdvig-v-K(G)} it follows immediately that $a_i\cdot a^{-1}\to 1$ in $\Spec{\mathcal K}(G)$. For the characteristic function $1_N\in{\mathcal K}(G)$ of the subgroup $N$ we have $1_N(a_i\cdot a^{-1})\to 1_N(1)=1$, hence starting from a certain index $i$ all elements $a_i\cdot a^{-1}$ lie in $N$. Take a compact set $S\subseteq {\mathcal K}(N)$. By Lemma \ref{LM:K(N)=K_N(G)} we can find a compact set $T\subseteq {\mathcal K}(G)$ (consisting of functions whose restrictions at $N$ lie in $S$), such that we obtain a bijection between $T$ and $S$. Since $a_i\cdot a^{-1}\to 1$ in $\Spec{\mathcal K}(G)$, we obtain
$$
v(a_i\cdot a^{-1})\underset{i\to\infty}{\underset{v\in T}{\rightrightarrows}}v(1)
$$
and this is equivalent to
$$
u(a_i\cdot a^{-1})\underset{i\to\infty}{\underset{u\in S}{\rightrightarrows}}u(1).
$$
This is true for each compact set $S\subseteq {\mathcal K}(N)$, hence $a_i\cdot a^{-1}\to 1$ in $\Spec{\mathcal K}(N)$. But in Lemma \ref{LM:Spec(R^n-times-K)} we already proved that $\Spec{\mathcal K}(N)=N$, hence we obtain that $a_i\cdot a^{-1}\to 1$ in $N$, and therefore, in $G$.
\epr

Theorem \ref{TH:Spec-K(G)=G} and \ref{C-obolochka-podalgebry-v-C(M)} immediately imply

\btm\label{TH:E(K(G))=C(G)}
If $G$ is a Moore group then the continuous envelope of the algebra ${\mathcal K}(G)$ is the algebra ${\mathcal C}(G)$:
\begin{equation}\label{E(K(G))=C(G)}
\Env_{\mathcal C} {\mathcal K}(G)={\mathcal C}(G)
\end{equation}
\etm

\subsection{$\Env_{\mathcal C}{\mathcal C}^\star(G)$ and ${\mathcal K}(G)$ as Hopf algebras}

\btm\label{TH:Env_C-C^*(G)-odot-Hopf}
If $G$ is a Moore group, then
\bit{
\item[(i)] the continuous envelope $\Env_{\mathcal C} {\mathcal C}^\star(G)$ of its group algebra ${\mathcal C}^\star(G)$ is an involutive Hopf algebra in the category of stereotype spaces $(\tt{Ste},\odot)$,

\item[(ii)] its dual algebra ${\mathcal K}(G)$ is an involutive Hopf algebra in the category of stereotype spaces $(\tt{Ste},\circledast)$.
}\eit
\etm
\bpr
We use here the fact that this proposition was proved in \cite{Akbarov-ker-coker} for the special case when $G$ is a compact buildup of an abelian group \cite[Theorem 2.8]{Akbarov-ker-coker}.

1. Consider first the case when $G$ is a Lie-Moore group. Then by theorem \ref{TH:Lie-Moore} $G$ is a finite extension of some compact buildup of an abelian (Lie) group:
$$
1\to Z\cdot K=N\to G\to F\to 1
$$
(here $Z$ is an abelian Lie group, $K$ a compact Lie group, $F$ a finite group). This chain generates a chain of homomorphisms of group algebras
$$
\C\to {\mathcal C}^\star(Z\cdot K)={\mathcal C}^\star(N)\to {\mathcal C}^\star(G)\to {\mathcal C}^\star(F)=\C_F\to \C.
$$
(the second equality in this chain follows from Example \ref{EX:Env_C-C_F=C_F}) and a chain of homomorphisms of smooth envelopes
$$
\C\to \Env_{\mathcal C} {\mathcal C}^\star(Z\cdot K)=\Env_{\mathcal C} {\mathcal C}^\star(N)\to \Env_{\mathcal C} {\mathcal C}^\star(G)\to \Env_{\mathcal C} {\mathcal C}^\star(F)=\eqref{Env_C-C_F=C_F}=\C_F\to \C.
$$

Consider the spaces ${\mathcal K}_L(G)$, $L\in G/N$, defined in \eqref{DEF:K_L(G)}. The sum of the characteristic functions $1_L$ is the identity of the algebra ${\mathcal K}(G)$,
$$
1=\sum_{L\in G/N}1_L,
$$
(here we have a finite number of summands). As a corollary, we can consider the space  ${\mathcal K}(G)$ as a finite sum of spaces ${\mathcal K}_L(G)$:
$$
{\mathcal K}(G)=\bigoplus_{L\in G/N}{\mathcal K}_L(G).
$$
Denote by $E_N$ the image of the space $\Env_{\mathcal C} {\mathcal C}^\star(N)$ in $\Env_{\mathcal C} {\mathcal C}^\star(G)$ under the mapping $\Env_{\mathcal C} {\mathcal C}^\star(N)\to \Env_{\mathcal C} {\mathcal C}^\star(G)$. Let for each $L\in G/N$ the symbol $E_L$ denote the shift of the space $E_N$ by any element $g_L\in L$:
$$
E_L=g_L*E_N.
$$
Since the shift is an automorphism of the space ${\mathcal K}(G)$, it is also an automorphism of the dual space $\Env_{\mathcal C} {\mathcal C}^\star(G)$, hence the space $E_L$ is well-defined. On the other hand, its definition does not depend on the choice of the element $g_L\in L$.

In these notations the dual space to each ${\mathcal K}_L(G)$ can be identified with the space $E_L$, since
$$
{\mathcal K}_L(G)^\star\cong {\mathcal K}(G)^\star/{\mathcal K}_L(G)^\perp\cong
\Env_{\mathcal C} {\mathcal C}^\star(G)/(\bigoplus_{M\ne L}E_M)\cong E_L.
$$
As a corollary, the space $\Env_{\mathcal C} {\mathcal C}^\star(G)$ can be represented as a direct sum of the spaces $E_L$:
$$
\Env_{\mathcal C}  {\mathcal C}^\star(G)=\bigoplus_{L\in G/N} E_L.
$$
Our aim is to show that this algebra is injective (i.e. is an algebra with respect to the tensor product $\odot$). First, note that by \cite[Theorem 2.8]{Akbarov-ker-coker}, the algebra $E_N=\Env_{\mathcal C}  {\mathcal C}^\star(N)$ is injective. Let
$$
m_N: E_N\odot E_N\to E_N
$$
be the continuous extension of its multiplication. We have to construct the operator
\begin{equation}\label{DEF:umnozhenie-na-Env-E*(G)-odot-Env-E*(G)}
m: \Env_{\mathcal C}  {\mathcal C}^\star(G)\odot \Env_{\mathcal C}  {\mathcal C}^\star(G)\to \Env_{\mathcal C}  {\mathcal C}^\star(G),
\end{equation}
which is an extension of the multiplication in $\Env_{\mathcal C}  {\mathcal C}^\star(G)$.

We can do this as follows. To each coset $L\in G/N$ we assign an arbitrary element $g_L\in L$. Consider the shift operators
$$
T^l_L:E_N\to E_L\quad\Big|\quad T^l_Lx=\env_{\mathcal C} \delta^{g_L}\cdot x,\qquad x\in E_N
$$
$$
T^r_L:E_N\to E_L\quad\Big|\quad T^r_Lx=x\cdot \env_{\mathcal C} \delta^{g_L},\qquad x\in E_N
$$
and the reverse shift operators
$$
(T^l_L)^{-1}:E_N\gets E_L\quad\Big|\quad (T^l_L)^{-1}y=\env_{\mathcal C} \delta^{g_L^{-1}}\cdot y,\qquad x\in E_L
$$
$$
(T^r_L)^{-1}:E_N\gets E_L\quad\Big|\quad (T^r_L)^{-1}y=y\cdot\env_{\mathcal C} \delta^{g_L^{-1}},\qquad x\in E_L.
$$
Let us decompose the injective tensor square of the space $\Env_{\mathcal C}  {\mathcal C}^\star(G)$ into the product of its components:
$$
\Env_{\mathcal C}  {\mathcal C}^\star(G)\odot \Env_{\mathcal C}  {\mathcal C}^\star(G)\cong
\bigoplus_{L,M\in G/N} E_L\odot E_M.
$$
Then the operator \eqref{DEF:umnozhenie-na-Env-E*(G)-odot-Env-E*(G)} can be defined on each component $E_L\odot E_M$ by the formula
\begin{equation}\label{PROOF:TH:gladk-obolochka-komp-porozhd-Lie-Moore-1}
m=T^l_L\circ T^r_M\circ m_N\circ \big((T^l_L)^{-1}\odot (T^r_M)^{-1}\big).
\end{equation}
To verify that the obtained operator extends the multiplication in $\Env_{\mathcal C}  {\mathcal C}^\star(G)$ it is sufficient to check this on delta-functionals (since their span is dense in ${\mathcal C}^\star(G)$ and therefore in $\Env_{\mathcal C}  {\mathcal C}^\star(G)$ as well). Take two elements $a,b\in G$ and find the cosets $L,M\in G/N$ such that $a\in L$ and $b\in M$. Then
\begin{multline*}
m(\env_{\mathcal C} \delta^a\odot\env_{\mathcal C} \delta^b)=\Big(T^l_L\circ T^r_M\circ m_N\circ \big((T^l_L)^{-1}\odot (T^r_M)^{-1}\big)\Big)(\env_{\mathcal C} \delta^a\odot\env_{\mathcal C} \delta^b)=\\=
\Big(T^l_L\circ T^r_M\circ m_N\Big)\big((T^l_L)^{-1}(\env_{\mathcal C} \delta^a)\odot (T^r_M)^{-1}(\env_{\mathcal C} \delta^b)\big)=\\=
\Big(T^l_L\circ T^r_M\circ m_N\Big)\big((\env_{\mathcal C} \delta^{g_L^{-1}}\cdot\env_{\mathcal C} \delta^a)\odot (\env_{\mathcal C} \delta^b\cdot\env_{\mathcal C} \delta^{g_M^{-1}})\big)=\\=
\Big(T^l_L\circ T^r_M\circ m_N\Big)\big(\env_{\mathcal C} (\delta^{g_L^{-1}}*\delta^a)\odot \env_{\mathcal C} (\delta^b*\delta^{g_M^{-1}})\big)=
\Big(T^l_L\circ T^r_M\Big)\big(\env_{\mathcal C} (\delta^{g_L^{-1}}*\delta^a)\cdot \env_{\mathcal C} (\delta^b*\delta^{g_M^{-1}})\big)=\\=
\Big(T^l_L\circ T^r_M\Big)\big(\env_{\mathcal C} (\delta^{g_L^{-1}}*\delta^a*\delta^b*\delta^{g_M^{-1}})\big)=
\env_{\mathcal C} (\delta^{g_L})\cdot\env_{\mathcal C} (\delta^{g_L^{-1}}*\delta^a*\delta^b*\delta^{g_M^{-1}})
\cdot\env_{\mathcal C} (\delta^{g_M})=\\=
\env_{\mathcal C} (\delta^{g_L}*\delta^{g_L^{-1}}*\delta^a*\delta^b*\delta^{g_M^{-1}}*\delta^{g_M})=
\env_{\mathcal C} (\delta^a*\delta^b)=
\env_{\mathcal C} (\delta^a)\cdot\env_{\mathcal C} (\delta^b)
\end{multline*}

2. Now let $G$ be an arbitrary Moore group. Then by Theorem \ref{TH:Moore=lim-Lie-Moore} it can be represented as a projective limit of a system of Lie-Moore groups:
$$
G=\leftlim_{\infty\gets i} G_i
$$
Each projection $G\to G_i$ induces a homomorphism of algebras of continuous functions ${\mathcal C}(G)\gets{\mathcal C}(G_i)$, and this system of homomorphisms generates the homomorphism of the injective limit
$$
{\mathcal C}(G)\gets\rightlim_{i\to\infty}{\mathcal C}(G_i).
$$
When we pass to the dual spaces we obtain a homomorphism into the projective limit
\begin{equation}\label{PROOF:nepr-obolochka-Moore-1}
{\mathcal C}^\star(G)\to\leftlim_{\infty\gets i}{\mathcal C}^\star(G_i).
\end{equation}
We have to show that this homomorphism generates the equality of the continuous envelopes:
\begin{equation}\label{PROOF:nepr-obolochka-Moore-2}
\Env_{\mathcal C} {\mathcal C}^\star(G)=\leftlim_{\infty\gets i}
\Env_{\mathcal C} {\mathcal C}^\star (G_i).
\end{equation}
Since $G_i$ are Lie-Moore groups, as we already proved, the algebras $\Env_{\mathcal C} {\mathcal C}^\star (G_i)$ are involutive Hopf algebras in the category $(\tt{Ste},\odot)$. By Corollary \ref{COR:proj-lim-Hopf-algebr-odot} we have that their projective limit $\Env_{\mathcal C} {\mathcal C}^\star (G)$ is an involutive Hopf algebra in the category $(\tt{Ste},\odot)$.

3. To prove \eqref{PROOF:nepr-obolochka-Moore-2}, let us first verify that the projective limit on the right is the continuous extension of the algebra ${\mathcal C}^\star(G)$. Consider arbitrary morphism $\ph: {\mathcal C}^\star(G)\to B$ into an arbitrary $C^*$-algebra $B$. The composition with the delta-functionals $\pi=\ph\circ\delta_G$ is a homomorphism of the Moore group into a Banach algebra $B$, hence it is factored through some projection $\rho_j:G\to G_i$:
$$
\xymatrix @R=2.pc @C=3.0pc 
{
G\ar@{-->}[r]^{\rho_j}\ar[dr]_{\pi} &  G_j\ar@{-->}[d]^{\pi_j} \\
 & B
}
$$
For each index $i\ge j$ we have the same:
$$
\xymatrix @R=2.pc @C=3.0pc 
{
G\ar@{-->}[r]^{\rho_i}\ar[dr]_{\pi} &  G_i\ar@{-->}[d]^{\pi_i} \\
 & B
}
$$
This gives a diagram for the algebras of measures:
$$
\xymatrix @R=2.pc @C=3.0pc 
{
{\mathcal C}^\star(G)\ar@{-->}[r]^{\sigma_i}\ar[dr]_{\ph} & {\mathcal C}^\star(G_i)\ar@{-->}[d]^{\ph_i} \\
 & B
}
$$
which can be extended to a diagram
$$
\xymatrix @R=2.pc @C=3.0pc 
{
{\mathcal C}^\star(G)\ar@{-->}[r]^{\sigma_i}\ar[dr]_{\ph} & {\mathcal C}^\star(G_i)\ar@{-->}[d]^{\ph_i}\ar@{-->}[r]^{\env_{\mathcal C} \sigma_i} &  \Env_{\mathcal C} {\mathcal C}^\star(G_i)\ar@{-->}[dl]^{\ph_i'} \\
 & B &
}
$$
If we throw out the vertex ${\mathcal C}^\star(G_i)$ and consider the product with $i\ge j$ we obtain the diagram
$$
\xymatrix @R=2.pc @C=3.0pc 
{
{\mathcal C}^\star(G)\ar@{-->}[rr]^{\sigma}\ar[dr]_{\ph} &  & \prod_{i\ge j} \Env_{\mathcal C} {\mathcal C}^\star(G_i)\ar@{-->}[dl]^{\prod_{i\ge j}\ph_i'} \\
 & B &
}
$$
Certainly, the image of the algebra ${\mathcal C}^\star(G)$ under the mapping $\sigma$ lies in the projective limit $\leftlim_{\infty\gets i, i\ge j} \Env_{\mathcal C} {\mathcal C}^\star (G_i)$, and we can replace the last diagram by the diagram
$$
\xymatrix @R=2.pc @C=3.0pc 
{
{\mathcal C}^\star(G)\ar@{-->}[rr]^{\sigma}\ar[dr]_{\ph} &  & \leftlim_{\infty\gets i, i\ge j} \Env_{\mathcal C} {\mathcal C}^\star (G_i) \ar@{-->}[dl]^{\leftlim_{\infty\gets i, i\ge j} \ph_i'} \\
 & B &
}
$$
After that we can notice that the limits $\leftlim_{\infty\gets i, i\ge j} \Env_{\mathcal C} {\mathcal C}^\star (G_i)$ and $\leftlim_{\infty\gets i} \Env_{\mathcal C} {\mathcal C}^\star (G_i)$ coincide, hence we obtain the diagram
$$
\xymatrix @R=2.pc @C=3.0pc 
{
{\mathcal C}^\star(G)\ar@{-->}[rr]^{\sigma}\ar[dr]_{\ph} &  & \leftlim_{\infty\gets i} \Env_{\mathcal C} {\mathcal C}^\star (G_i) \ar@{-->}[dl] \\
 & B &
}
$$

4. We understood that the projective limit on the right in \eqref{PROOF:nepr-obolochka-Moore-2} is a continuous extension of the algebra ${\mathcal C}^\star(G)$. Now we have to verify that it is a continuous envelope. Let  $\varepsilon:{\mathcal C}^\star(G)\to \Env_{\mathcal C}{\mathcal C}^\star(G)$ be another continuous extension.

Consider a Lie-Moore group $G_i$ and the corresponding projection $\sigma_i:{\mathcal C}^\star(G)\to {\mathcal C}^\star(G_i)$. Take an arbitrary $C^*$-neighborhood of zero $U$ in ${\mathcal C}^\star(G_i)$, and let $\pi^U:{\mathcal C}^\star(G_i)\to {\mathcal C}^\star(G_i)/U$ be the corresponding $C^*$-quotient mapping (we use the terminology of \cite{Akbarov-env} and \cite{Akbarov-C^infty-2}). Since $\env_{\mathcal C}{\mathcal C}^\star(G):{\mathcal C}^\star(G)\to \Env_{\mathcal C}{\mathcal C}^\star(G)$ is a continuous extension, in the diagram
$$
\xymatrix @R=2.pc @C=3.0pc 
{
& {\mathcal C}^\star(G)\ar[ld]_{\env_{\mathcal C}{\mathcal C}^\star(G)}\ar[rd]^{\sigma_i} & \\
\Env_{\mathcal C}{\mathcal C}^\star(G)\ar@{-->}[rd]_{\ph_i^U} &  & {\mathcal C}^\star (G_i)\ar[ld]^{\pi^U}  \\
& {\mathcal C}^\star (G_i)/U &
}
$$
there is a unique dashed arrow $\ph_i^U$. This is true for each $C^*$-neighborhood of zero $U$ in ${\mathcal C}^\star(G_i)$, and it is easy to see, when we make $U$ smaller, the corresponding arrows $\ph_i^U$ are connected to each other with natural mediators (here  $V\subseteq U$):
$$
\xymatrix @R=2.pc @C=3.0pc 
{
& {\mathcal C}^\star(G)\ar[ld]_{\env_{\mathcal C}{\mathcal C}^\star(G)}\ar[rd]^{\sigma_i} & \\
\Env_{\mathcal C}{\mathcal C}^\star(G)\ar@{-->}@/_5ex/[rdd]_{\ph_i^U}\ar@{-->}[rd]_{\ph_i^V} &  & {\mathcal C}^\star (G_i)\ar[ld]^{\pi^V}\ar@/^5ex/[ldd]^{\pi^U}  \\
& {\mathcal C}^\star(G_i)/V\ar@{-->}[d] & \\
& {\mathcal C}^\star(G_i)/U &
}
$$
This means that we can pass to the projective limit, and we obtain the diagram
$$
\xymatrix @R=2.pc @C=3.0pc 
{
& {\mathcal C}^\star(G)\ar[ld]_{\env_{\mathcal C}{\mathcal C}^\star(G)}\ar[rd]^{\sigma_i} & \\
\Env_{\mathcal C}{\mathcal C}^\star(G)\ar@{-->}[rd]_{\ph_i} &  & {\mathcal C}^\star (G_i)\ar[ld]^{\pi}  \\
& \lim\limits_{0\gets U}{\mathcal C}^\star (G_i)/U &
}
$$
Here $\env_{\mathcal C}{\mathcal C}^\star(G)$ is a dense epimorphism, i.e. an epimorphism in the category of stereotype spaces. By Lemma \ref{LM:epimorphism} this diagram can be complemented to the diagram
$$
\xymatrix @R=2.pc @C=3.0pc 
{
& {\mathcal C}^\star(G)\ar[ld]_{\env_{\mathcal C}{\mathcal C}^\star(G)}\ar[rd]^{\sigma_i}\ar@{-->}[d] & \\
\Env_{\mathcal C}{\mathcal C}^\star(G)\ar[rd]_{\ph_i}\ar@{-->}[r] & \Im_\infty(\pi\circ\sigma_i)\ar@{-->}[d] & {\mathcal C}^\star (G_i)\ar[ld]^{\pi}  \\
& \lim\limits_{0\gets U}{\mathcal C}^\star (G_i)/U &
}
$$
Recall now that $\sigma_i$ is also a dense epimorphism. Hence by \eqref{Im_infty-ph-circ-e=Im_infty-ph}
$$
\Im_\infty(\pi\circ\sigma_i)=\Im_\infty\pi.
$$
Let us put this into the diagram:
$$
\xymatrix @R=2.pc @C=3.0pc 
{
& {\mathcal C}^\star(G)\ar[ld]_{\env_{\mathcal C}{\mathcal C}^\star(G)}\ar[rd]^{\sigma_i}\ar[d] & \\
\Env_{\mathcal C}{\mathcal C}^\star(G)\ar[rd]_{\ph_i}\ar[r] & \Im_\infty\pi\ar[d] & {\mathcal C}^\star (G_i)\ar[ld]^{\pi}\ar@{-->}[l]  \\
& \lim\limits_{0\gets U}{\mathcal C}^\star (G_i)/U &
}
$$
Now by formula \cite[(5.61)]{Akbarov-env} we have
$$
\Im_\infty\pi=\Env_{\mathcal C} {\mathcal C}^\star(G_i).
$$
Put this into the diagram:
$$
\xymatrix @R=2.pc @C=3.0pc 
{
& {\mathcal C}^\star(G)\ar[ld]_{\env_{\mathcal C}{\mathcal C}^\star(G)}\ar[rd]^{\sigma_i}\ar[d] & \\
\Env_{\mathcal C}{\mathcal C}^\star(G)\ar[rd]_{\ph_i}\ar[r] & \Env_{\mathcal C} {\mathcal C}^\star(G_i)\ar[d] & {\mathcal C}^\star (G_i)\ar[ld]^{\pi}\ar@{-->}[l]_{\env_{\mathcal C} {\mathcal C}^\star(G_i)}  \\
& \lim\limits_{0\gets U}{\mathcal C}^\star (G_i)/U &
}
$$
If we throw away the lower and the right nodes, we obtain
$$
\xymatrix @R=2.pc @C=3.0pc 
{
& {\mathcal C}^\star(G)\ar[ld]_{\env_{\mathcal C}{\mathcal C}^\star(G)}\ar[d]^{\env_{\mathcal C} {\mathcal C}^\star(G_i)\circ\sigma_i}  \\
\Env_{\mathcal C}{\mathcal C}^\star(G)\ar[r] & \Env_{\mathcal C} {\mathcal C}^\star(G_i)  \\
}
$$
Now let us change the index $i$. For $i\le j$ we have the diagram
$$
\xymatrix @R=2.pc @C=3.0pc 
{
& {\mathcal C}^\star(G)\ar[ld]_{\env_{\mathcal C}{\mathcal C}^\star(G)}\ar[d]^{\env_{\mathcal C} {\mathcal C}^\star(G_j)\circ\sigma_j}\ar@/^8ex/[ddr]^{\env_{\mathcal C} {\mathcal C}^\star(G_i)\circ\sigma_i} &  \\
\Env_{\mathcal C}{\mathcal C}^\star(G)\ar[r]\ar@/_5ex/[rrd] & \Env_{\mathcal C} {\mathcal C}^\star(G_j)\ar[dr] & \\
& & \Env_{\mathcal C} {\mathcal C}^\star(G_i)
}
$$
As a corollary, we can pass to the projective limit by indices $i$:
$$
\xymatrix @R=2.pc @C=3.0pc 
{
& {\mathcal C}^\star(G)\ar[ld]_{\env_{\mathcal C}{\mathcal C}^\star(G)}\ar[d]^{\lim\limits_{\infty\gets i}\env_{\mathcal C} {\mathcal C}^\star(G_i)\circ\sigma_i}  \\
\Env_{\mathcal C}{\mathcal C}^\star(G)\ar@{-->}[r] & \lim\limits_{\infty\gets i}\Env_{\mathcal C} {\mathcal C}^\star(G_i)  \\
}
$$
The dashed arrow in this diagram is the very same arrow which shows that $\lim\limits_{\infty\gets i}\Env_{\mathcal C} {\mathcal C}^\star(G_i)$ is a continuous envelope.
\epr

\section{Continuous duality}

\subsection{Reflexivity with respect to the envelope.}

We say that an involutive stereotype Hopf algebra $H$  in the category $({\tt Ste},\circledast)$ is {\it continuously reflexive}\label{DEF:reflexiv-otn-obolochki}, if its continuous envelope $\Env_{\mathcal C} H$ has a structure of injective Hopf algebra in the category $({\tt Ste},\odot)$ such that the following two requirements hold:
\bit{
\item[(i)] a morphism of the continuous envelope $\env_{\mathcal C} H:H\to \Env_{\mathcal C} H$ is a homomorphism of Hopf algebras in the sense that the following diagrams are commutative:
\begin{equation}\label{DIAG:reflex-otn-obolochki-1}
 \xymatrix @R=2.pc @C=2.pc
{
& H\odot H\ar[dr]^{\quad \env_{\mathcal C} H\odot \env_{\mathcal C} H} & \\
H\circledast H\ar[ur]^{@}\ar[dr]^{\quad \env_{\mathcal C} H\circledast \env_{\mathcal C} H}\ar[dd]_{\mu} & & \Env_{\mathcal C} H\odot \Env_{\mathcal C} H\ar[dd]_{\mu_E} \\
& \Env_{\mathcal C} H\circledast \Env_{\mathcal C} H\ar[ur]^{@} & \\
H\ar[rr]^{\env_{\mathcal C} H} && \Env_{\mathcal C} H
}
\end{equation}
\begin{equation}\label{DIAG:reflex-otn-obolochki-2}
 \xymatrix @R=2.pc @C=2.pc
{
& H\odot H\ar[dr]^{\quad\env_{\mathcal C} H\odot \env_{\mathcal C} H} & \\
H\circledast H\ar[ur]^{@}\ar[dr]^{\quad \env_{\mathcal C} H\circledast \env_{\mathcal C} H} & & \Env_{\mathcal C} H\odot \Env_{\mathcal C} H \\
& \Env_{\mathcal C} H\circledast \Env_{\mathcal C} H\ar[ur]^{@} & \\
H\ar[rr]^{\env_{\mathcal C} H}\ar[uu]^{\varkappa} && \Env_{\mathcal C} H\ar[uu]^{\varkappa_E}
}
\end{equation}
\begin{equation}\label{DIAG:reflex-otn-obolochki-3}
 \xymatrix @R=2.pc @C=2.pc
{
H\ar[rr]^{\env_{\mathcal C} H} & & \Env_{\mathcal C} H\\
& \C\ar[ul]^{\iota}\ar[ur]_{\iota_E} &
}\qquad
 \xymatrix @R=2.pc @C=2.pc
{
H\ar[rr]^{\env_{\mathcal C} H}\ar[dr]_{\e} & & \Env_{\mathcal C} H\ar[dl]^{\e_E} \\
& \C &
}
\end{equation}
\begin{equation}\label{DIAG:reflex-otn-obolochki-4}
 \xymatrix @R=3.pc @C=4.pc
{
H\ar[r]^{\env_{\mathcal C} H}\ar[d]_{\sigma} & \Env_{\mathcal C} H\ar[d]^{\sigma_E} \\
H\ar[r]^{\env_{\mathcal C} H} & \Env_{\mathcal C} H
}\qquad
 \xymatrix @R=3.pc @C=4.pc
{
H\ar[r]^{\env_{\mathcal C} H}\ar[d]_{\bullet} & \Env_{\mathcal C} H\ar[d]^{\bullet_E} \\
H\ar[r]^{\env_{\mathcal C} H} & \Env_{\mathcal C} H
}
\end{equation}
-- here $@$ is the Grothendieck transformation\footnote{See \cite{Akbarov} or \cite{Akbarov-C^infty-1}}, $\mu$, $\iota$, $\varkappa$, $\e$, $\sigma$, $\bullet$ are structural morphisms (multiplication, unity, comultiplication, counity, antipode, involution) in $H$, and $\mu_E$, $\iota_E$, $\varkappa_E$, $\e_E$, $\sigma_E$, $\bullet_E$ the structural morphisms in $\Env_{\mathcal C} H$.

\item[(ii)]\label{(env-H)^star:H^star-gets-(Env-H)^star} the mapping $(\env_{\mathcal C} H)^\star:H^\star\gets (\Env_{\mathcal C} H)^\star$ dual to the morphism of continuous envelope $\env_{\mathcal C} H:H\to\Env_{\mathcal C} H$, is again a continuous envelope:
$$
(\env_{\mathcal C} H)^\star=\env_{\mathcal C} (\Env_{\mathcal C} H)^\star
$$

}\eit

It is convenient to display the conditions (i) and (ii) as a diagram
 \begin{equation}\label{obshaya-diagramma-refleksivnosti}
 \xymatrix @R=1.pc @C=1.pc
 {
 H
 & \ar@{|->}[r]^{\env_{\mathcal C}} & &
\Env_{\mathcal C} H
 \\
 & & &
 \ar@{|->}[d]^{\star}
 \\
 \ar@{|->}[u]^{\star}
 & & &
 \\
 H^\star
 & &
 \ar@{|->}[l]_{\env_{\mathcal C}}
 &
 (\Env_{\mathcal C} H)^\star
 }
 \end{equation}
which we call the {\it reflexivity diagram}, and which we endow with the following sense:
 \bit{
\item[1)] in the corners of the square there are involutive Hopf algebras; the first algebra, $H$, is the Hopf algebra in $({\tt Ste},\circledast)$, the second algebra, $\Env_{\mathcal C} H$, is the Hopf algebra in $({\tt Ste},\odot)$, and further the categories $({\tt Ste},\circledast)$ and $({\tt Ste},\odot)$ alternate,

\item[2)] the alternation of the operations $\env_{\mathcal C}$ and $\star$ (no matter where we start) on the fourth step returns us back to the initial Hopf algebra (certainly, up to an isomorphism of functors).
 }\eit

The sense of the term ``reflexivity'' here is as follows. Denote the single successive application of the operations $\env$ and $\star$ by some symbol, for example, $\widehat{\ }\,\,$,
$$
\widehat H:=(\Env_{\mathcal C} H)^\star
$$
Since $\Env_{\mathcal C} H$ has a unique structure of Hopf algebra with respect to $\odot$, the dual space $\widehat H=(\Env_{\mathcal C} H)^\star$ has a structure of involutive Hopf algebra with respect to $\circledast$. Moreover, $\widehat H=(\Env_{\mathcal C} H)^\star$ is a Hopf algebra, reflexive with respect to $\Env_{\mathcal C}$, since the application of $\star$ to the diagrams \eqref{DIAG:reflex-otn-obolochki-1}-\eqref{DIAG:reflex-otn-obolochki-4} gives the same diagrams with the replacement $H$ by $\widehat H=(\Env_{\mathcal C} H)^\star$ (we use here the condition (ii) on page \pageref{(env-H)^star:H^star-gets-(Env-H)^star}).

Let us call $\widehat H=(\Env_{\mathcal C} H)^\star$ the {\it dual Hopf algebra to
$H$ with respect to the envelope $\Env_{\mathcal C}$}. The diagram \eqref{obshaya-diagramma-refleksivnosti} means that $H$ is naturally isomorphic to its second dual Hopf algebra in this sense:
 \begin{equation}\label{H-cong-(H^*)^*}
H\cong \widehat{\widehat H}
 \end{equation}

\subsection{Continuous reflexivity.} Theorems \ref{TH:E(K(G))=C(G)} and \ref{TH:Env_C-C^*(G)-odot-Hopf} imply the following main result of our work (and in the corrected formulations, the main result of the Yu.~N.~Kuznetsova work \cite{Kuznetsova}):

 \btm\label{TH:nepr-dvoistvennost}
If $G$ is a Moore group, then the algebras ${\mathcal C}^\star(G)$ and ${\mathcal
K}(G)$ are continuously reflexive, and the reflexivity diagram for them is  \eqref{chetyrehugolnik-C-C*-0}:
 \begin{equation}\label{chetyrehugolnik-C-C*}
 \xymatrix @R=1.pc @C=2.pc
 {
 {\mathcal C}^\star(G)
 & \ar@{|->}[r]^{\Env_{\mathcal C} } & &
 \Env_{\mathcal C} {\mathcal C}^\star(G)
 \\
 & & &
 \ar@{|->}[d]^{\star}
 \\
 \ar@{|->}[u]^{\star}
 & & &
 \\
 {\mathcal C}(G)
 & &
 \ar@{|->}[l]_{\Env_{\mathcal C} }
 &
 {\mathcal K}(G)
 }
 \end{equation}
 \etm
\bpr
When we move by the chain \eqref{chetyrehugolnik-C-C*} from the left lower corner, ${\mathcal C}(G)$, we come to the algebra ${\mathcal K}(G)$, which by Theorem \ref{TH:E(K(G))=C(G)} turns into the algebra ${\mathcal C}(G)$ under the action of the envelope $\Env_{\mathcal C}$, and thus the chain \eqref{chetyrehugolnik-C-C*} closes. On the other hand, in this diagram the elements ${\mathcal C}(G)$ and $\Env_{\mathcal C}{\mathcal C}^\star(G)$ are $\odot$-Hopf algebras (${\mathcal C}(G)$ due to \cite[Example 10.24]{Akbarov} and \cite[4.2]{Akbarov-stein-groups}, and $\Env_{\mathcal C}{\mathcal C}^\star(G)$ by Theorem \ref{TH:Env_C-C^*(G)-odot-Hopf}), and the elements ${\mathcal C}^\star(G)$ and ${\mathcal K}(G)$ are $\circledast$-Hopf algebras (${\mathcal C}^\star(G)$ due to \cite[Example 10.24]{Akbarov} and \cite[4.2]{Akbarov-stein-groups}, and ${\mathcal K}(G)$ again by Theorem \ref{TH:Env_C-C^*(G)-odot-Hopf}).
\epr

\tableofcontents


\begin{thebibliography}{99}



\bibitem{Akbarov} S.~S.~Akbarov. Pontryagin duality in the theory of topological vector spaces and in topological algebra, Journal of Mathematical Sciences, 113(2):179-349, 2003.


\bibitem{Akbarov-stein-groups} S.~S.~Akbarov. Holomorphic functions of exponential type and duality for Stein groups with algebraic connected component of identity, {\it
Journal of Mathematical Sciences}, 162(4): 459-586, 2009;
\href{http://arxiv.org/abs/0806.3205}{http://arxiv.org/abs/0806.3205}.

\bibitem{Akbarov-env} S.~S.~Akbarov. Envelopes and refinements in categories, with applications to functional analysis. Dissertaciones mathematicae, 513(1): 1-188, 2016, \href{https://www.impan.pl/en/publishing-house/journals-and-series/dissertationes-mathematicae/all/513}{https://www.impan.pl/en/publishing-house/journals-and-series/dissertationes-mathematicae/all/513}, \href{http://arxiv.org/abs/1110.2013}{http://arxiv.org/abs/1110.2013}.


\bibitem{Akbarov-C^infty-1} S.~S.~Akbarov. Continuous and smooth envelopes of topological algebras. Part I; Journal of Mathematical Sciences, 227(5):531-668, 2017; \href{https://arxiv.org/abs/1303.2424}{https://arxiv.org/abs/1303.2424}.


\bibitem{Akbarov-C^infty-2} S.~S.~Akbarov. Continuous and smooth envelopes of topological algebras. Part II; Journal of Mathematical Sciences, 227(6):669-789, 2017; \href{https://arxiv.org/abs/1303.2424}{https://arxiv.org/abs/1303.2424}.


\bibitem{Akbarov-ker-coker} S.~S.~Akbarov. Kernel and cokernel in the category of augmented stereotype algebras.


\bibitem{Pirkovskii-DM} A.~Yu.~Pirkovskii. Stably flat completions of universal enveloping algebras. Dissertationes Math. (Rozprawy Math.) 441 (2006), 1--60.

\bibitem{Pirkovskii-PAMS} A.~Yu.~Pirkovskii. Arens-Michael enveloping algebras and analytic smash products. Proc. Amer. Math. Soc. 134 (2006), no.9, 2621--2631.

\bibitem{Pirkovskii-MMO} A.~Yu.~Pirkovskii. Arens–Michale envelopes, homological epimorphisms, and relatively quasi-free algebras. Trans. Moscow Math. Soc. 69: 27–104, 2008.

\bibitem{Grosser-Moskowitz-2} S.~Grosser, M.~Moskowitz, Compactness conditions in topological groups. J. Reine Angew. Math., 246: 1–40, 1971.

\bibitem{Kuznetsova} Yu.~Kuznetsova, A duality for Moore groups. J. Oper. Theory, 69(2):101-130, 2013, \href{http://arxiv.org/abs/0907.1409}{http://arxiv.org/abs/0907.1409}.


\bibitem{Morris} S.~A.~Morris, Pontryagin Duality and the Structure of Locally Compact Abelian Groups. Cambridge University Press, 1977.


\bibitem{Palmer} Th.~W.~Palmer. Banach algebras and the general theory of *-algebras. Vol. II. Academic Press. 2001.




\end{thebibliography}
\end{document}